\documentclass[a4j,10pt,leqno]{amsart}
\usepackage{amsmath,amssymb,amscd, amsthm}
\usepackage{epic,eepic,epsfig}
\usepackage{comment} 
\usepackage{mathrsfs}
\usepackage[all,cmtip]{xy}
\usepackage{graphicx}
\usepackage{here}
\usepackage{tikz}
\usepackage{tikz-cd}
\usepackage{pgfplots}
\usepackage{mathtools}
 \usepackage[nodayofweek]{datetime}
\makeatletter
\def\@settitle{\begin{center}%
  \baselineskip14\p@\relax
  \normalfont\LARGE\bfseries
  \@title
  \ifx\@subtitle\@empty\else
     \\[1ex] %\uppercasenonmath\@subtitle
     \normalsize\mdseries\@subtitle
  \fi
 \ifx\@didication\@empty\else
     \\[2ex] %\uppercasenonmath\@subtitle
     \large\mdseries\it\@dedication
  \fi
  \end{center}%
}
\def\subtitle#1{\gdef\@subtitle{#1}}
\def\@subtitle{}
\def\dedication#1{\gdef\@dedication{#1}}
\def\@dedication{}
\makeatother

\usepackage{wrapfig, url}
\usepackage[%
 setpagesize=false,%
 bookmarks=true,%
 bookmarksdepth=tocdepth,%
 bookmarksnumbered=true,%
 colorlinks=false,%
 pdftitle={},%
 pdfsubject={},%
 pdfauthor={},%
 pdfkeywords={}%
]{hyperref}

\makeatletter

\makeatletter
\renewcommand{\section}{\@startsection
{section}{1}{0mm}{5mm}{2mm}{\raggedright\bfseries}}
\makeatother
\newtheorem{theorem}{Theorem}[section] %Theorem (1.1)

\theoremstyle{definition}

\newtheorem{nsTheorem}[theorem]{Theorem} %{\bf}{\it}
\newtheorem{nsLemma}[theorem]{Lemma} %{\bf}{\it}
\newtheorem{nsProposition}[theorem]{Proposition} %{\bf}{\it}
\newtheorem{nsDefinition}[theorem]{Definition} %{\bf}{\rm}
\newtheorem{nsRemark}[theorem]{Remark} %{\bf}{\rm}
\newtheorem{nsNotation}[theorem]{Notation} %{\bf}{\rm}
\allowdisplaybreaks[1]
\begin{document}

\title{%[draft] 
Landen's trilogarithm functional equation \\
and 
$\ell$-adic Galois multiple polylogarithms}
%\subtitle{Here can be put a subtitile}  %option
\dedication{In memory of Toshie Takata}  %option

\author{Hiroaki Nakamura }


\address{Hiroaki Nakamura: 
Department of Mathematics, 
Graduate School of Science, 
Osaka University, 
Toyonaka, Osaka 560-0043, Japan}
\email{nakamura@math.sci.osaka-u.ac.jp}

\author{Densuke Shiraishi}
\address{Densuke Shiraishi: 
Department of Mathematics,
Faculty of Science Division II,
Tokyo University of Science,
1-3 Kagurazaka,
Shinjuku,
Tokyo 162-8601,
Japan}
\email{dshiraishi@rs.tus.ac.jp,
densuke.shiraishi@gmail.com}

\maketitle

%%%%%%%%%%%%%%%%%%%%%%%%%%%%%%%%%%%%%%%%%%%%%%%%%%%%%%%%%%%%%%%%%%%%
% This defines a short-running title of this paper !! 
\markboth{H.Nakamura, D.Shiraishi}
{Landen's trilogarithm functional equation}
%%%%%%%%%%%%%%%%%%%%%%%%%%%%%%%%%%%%%%%%%%%%%%%%%%%%%%%%%%%%%%%%%%

\subtitle{\it\large \hfill
In memory of Toshie Takata
} 
\maketitle

\begin{abstract}
{The Galois action on the pro-$\ell$ \'etale fundamental groupoid 
of the projective line minus three points with rational base points gives rise to
 a non-commutative formal power series in two variables with $\ell$-adic coefficients,
called the $\ell$-adic Galois associator.
In the present paper,
we focus on how Landen's functional equation of trilogarithms
and its $\ell$-adic Galois analog can be derived algebraically from  
the $S_3$-symmetry of the projective line minus three points.
Twofold proofs of the functional equation will be presented, 
one is based on the chain rule for the associator power series and the other 
is based on Zagier's tensor criterion devised in the framework of graded Lie algebras.
In the course of the second proof, we are led to investigate 
$\ell$-adic Galois multiple polylogarithms appearing
as regular coefficients of the $\ell$-adic Galois associator.
As an application, we show an $\ell$-adic Galois analog of 
Oi-Ueno's functional equation between 
$Li_{1,\dots,1,2}(1-z)$ and $Li_k(z)$'s $(k=1,2,...)$ .
}
\end{abstract}

\tableofcontents

\footnote[0]{
This paper has been published with open accessible format at 
\url{https://doi.org/10.1007/978-981-97-3778-9_8} 
\\
in 
``Low Dimensional Topology and Number Theory'' 
(M.Morishita, H.Nakamura, J.Ueki eds.), 
Springer Proceedings in Mathematics \& Statistics (PROMS, volume 456), pages 237--262, Springer Nature Singapore, 2025. 
% with open accessible format at 
%\url{https://doi.org/10.1007/978-981-97-3778-9_8}. 
%\quad\LaTeX -compiled on: \today \ at \currenttime 
}

\section{Introduction}

The study of polylogarithms, especially their functional equations, 
originated in the late $18$th century by Euler, Landen, and others.
The classical polylogarithm they studied is 
a complex function defined by the following power series 
\[
Li_k(z):=\frac{z}{1^k}+\frac{z^2}{2^k}+\frac{z^3}{3^k}+\cdots~~(|z|<1).
\]
For $k=2$,
it is called the dilogarithm,
and for $k=3$,
it is called the trilogarithm.
The multiple polylogarithm $Li_{\mathbf{k}}(z)$ 
for a multi-index $\mathbf{k}=(k_1\dots,k_d)\in\mathbb{N}^d$
generalizes $Li_k(z)$, which is defined by the power series
\[
Li_{\mathbf{k}}(z):=\sum_{0<n_1<\cdots<n_d}\frac{z^{n_d}}{n_1^{k_1}\cdots n_d^{k_d}}\quad (|z|<1).
\]
Note that $Li_k(z)=Li_{(k)}(z)$.
The collection $\{Li_{\mathbf{k}}(z)\}_{\mathbf{k}}$ satisfies certain recursive differential equations, from which follows that each
$Li_{\mathbf{k}}(z)$ has an iterated integral
expression that can be analytically continued 
to the universal cover 
of the three punctured Riemann sphere 
$\mathbf{P}^1(\mathbb{C})-\{0,1,\infty\}$
(cf. e.g. \cite[Lemma 1.5]{F04}).
There are known a number of functional equations between these functions 
evaluated at points with 
suitably chosen tracking paths from the unit segment $(0,1)$ on 
$\mathbf{P}^1(\mathbb{C}) - \{0,1,\infty\}$.
For example, the following formulas are typical:
\begin{align} \label{EulerLi2}
&Li_2(z)+Li_2(1-z) = \zeta(2)-\log(z)\log(1-z) ,\\
%\quad [{\rm Euler}, 1768], 
&Li_2(z)+Li_2\left(\frac{z}{z-1}\right)=-\frac12\log^2(1-z), 
\label{LandenLi2}
%\tag{$\#_{\mathbb{C}.{\rm Landen}}$} \\
\\
&Li_3(z)+Li_3(1-z)+Li_3\left(\frac{z}{z-1}\right)
\label{LandenLi3} \\
&\qquad = \zeta(3)+\zeta(2)\log(1-z)-\frac{1}{2}\log(x)\log^{2}(1-z)
 +\frac{1}{6}\log^{3}(1-z) .
\notag
%\quad [{\rm Landen}, 1780].
\end{align}
The former (\ref{EulerLi2}) is due to Leonhard Euler 
\cite{Euler1768} and the latter two (\ref{LandenLi2})-(\ref{LandenLi3})
are due to John Landen \cite{Landen1780}. 
See Lewin's book \cite{L81} for many other functional equations for polylogarithms.
As for multiple polylogarithms, in \cite{Oi2009}-\cite{Oi-Ueno2013}, Shu Oi and Kimio Ueno showed the following 
functional equation:
\begin{equation}
\label{OiUeno}
\sum_{j=0}^{k-1}Li_{k-j}(z)\frac{(-\log z)^j}{j!}+
Li_{\hspace{-2mm}\underbrace{\scriptstyle 1,\ldots,1}_{k-2~\text{times}}\hspace{-2mm},2}~(1-z)=\zeta(k)
\qquad (k\ge 2).
% \quad [{\rm Oi\text{-}Ueno}, 2013],[{\rm Oi}, 2009]
%\tag{$\#_{\mathbb{C}.1}$}.
\end{equation}

Let $\ell$ be a fixed prime.
The $\ell$-adic Galois multiple polylogarithm 
\[
{Li}^\ell_{\mathbf{k}}(z)
\left(={Li}^\ell_{\mathbf{k}}(\gamma_z:\overrightarrow{01}{\leadsto} z) 
\right): G_K \to \mathbb{Q}_{\ell}
\]
%${Li}^{\hspace{0.04cm}\ell,z}_{\mathbf{k}}$ 
is
a function on the absolute Galois group $G_K:={\rm Gal}(\overline{K}/K)$ 
of a subfield $K$ of $\mathbb{C}$ defined, for $\mathbf{k}=(k_1\dots,k_d)\in\mathbb{N}^d$ 
and an $\ell$-adic \'etale path $\gamma_z$ from $\overrightarrow{01}$ to a 
$K$-rational (tangential) point $z$ on $\mathbb{P}^1-\{0,1,\infty\}$, 
as a certain (signed) coefficient of the non-commutative formal power series
\begin{equation}
\label{Galois-associator-intro}
{\mathfrak f}_\sigma^{\gamma_z}(X,Y) \in \mathbb{Q}_{\ell}\langle\!\langle X,Y\rangle\!\rangle \quad (\sigma \in G_K)
\end{equation}
called the $\ell$-adic Galois associator. 
The functions ${Li}^\ell_{\mathbf{k}}(z)$ were originally introduced 
and called {\it the $\ell$-adic iterated integrals}
in a series of papers by Zdzis{\l}aw Wojtkowiak 
 (cf.\,e.g.,\,\cite{W0}-\cite{W3}).
In particular,
\begin{equation} \label{l-adic_zeta}
{\boldsymbol \zeta}_{{\mathbf k}}^\ell(\sigma):=
{Li}^\ell_{\mathbf{k}}(\delta:\overrightarrow{01}{\leadsto} {\overrightarrow{10}})(\sigma)
\end{equation}
for the standard path $\delta$ along the unit interval $(0,1)\subset{\mathbb R}$.
For $z\in K$ with a path $\gamma_z:\overrightarrow{01}{\leadsto} z$, we also write
\[
\rho_{z}(=\rho_{\gamma_z}): G_K \to {\mathbb Z}_{\ell}
\]
for the Kummer $1$-cocycle of the $\ell$-th power roots $\{z^{1/{\ell^n}}\}_{n}$
determined by $\gamma_z$.

In \cite{NW12}, Wojtkowiak and the first named author of the present paper
devised Zagier's tensor criterion for functional equations as a means
to calculate exact forms of identities with lower degree terms  
for both complex and $\ell$-adic Galois polylogarithms. 
Applying the method, we established a few examples of functional equations 
in both polylogarithms. In particular, the above (\ref{EulerLi2}) and (\ref{LandenLi2})
were shown to have the following $\ell$-adic Galois counterparts:
\begin{align}
Li_2^\ell(z)(\sigma)+Li_2^\ell(1-z)(\sigma) &=
{\boldsymbol \zeta}_{{2}}^\ell(\sigma) 
-\rho_z(\sigma)\rho_{1-z}(\sigma) ,
\label{GaloisLi2FE}
\\
Li_2^\ell(z)(\sigma)+Li_2^\ell \left(\frac{z}{z-1}\right)(\sigma) 
&=-\frac{\rho_{1-z}(\sigma)^2+\rho_{1-z}(\sigma)}{2}
\label{GaloisLi2Landen}
\end{align}
for $\sigma\in G_K$
(cf.\,\cite{NW12}.
See (\ref{GaloisTilchi2FE})-(\ref{GaloisTilchi2Landen})
and Proposition \ref{compareGaloisPolylog}
below for adjustments of notations.)

The purpose of this paper is to provide algebraic proofs of 
(\ref{LandenLi3}) and (\ref{OiUeno})
which can be used to obtain their $\ell$-adic Galois analogs 
reading as follows:

\begin{nsTheorem}[$\ell$-adic Galois analog of the Landen trilogarithm functional equation]
\label{l-adic Landen}
There are suitable paths $\overrightarrow{01}{\leadsto} 1-z, \overrightarrow{01}{\leadsto}\frac{z}{z-1}$ associated to
a given path $\gamma_z:\overrightarrow{01}{\leadsto} z$ such that 
the following functional equation
\begin{align*}
\label{GaloisLandenLi3}
&{Li}^\ell_{3}(z)(\sigma)+{Li}^\ell_{3}(1-z)(\sigma)
+{Li}^\ell_{3}\left(\frac{z}{z-1}\right)(\sigma) 
\\
%\tag{$\#_{\ell.{\rm Landen}}$}
&={\boldsymbol \zeta}_{{3}}^\ell(\sigma)-{\boldsymbol \zeta}_{{2}}^\ell(\sigma)\rho_{1-z}(\sigma)+\frac{1}{2}\rho_{z}(\sigma){\rho_{1-z}(\sigma)}^2-\frac{1}{6}{\rho_{1-z}(\sigma)}^3 
\\
&\qquad
-\frac{1}{2}{Li}^\ell_{2}(z)(\sigma)
-\frac{1}{12}\rho_{1-z}(\sigma)-\frac{1}{4}{\rho_{1-z}(\sigma)}^2
\end{align*}
holds for $\sigma\in G_K$.
\end{nsTheorem}

\begin{nsTheorem}[$\ell$-adic Galois analog of the Oi-Ueno functional equation]
\label{GaloisOiUeno}
\[ 
\sum_{j=0}^{k-1}{Li}^\ell_{k-j}(z)(\sigma)\frac{\rho_z(\sigma)^j}{j!} 
+{Li}^\ell_{\hspace{-2mm}\underbrace{\scriptstyle 1,\ldots,1}_{k-2~\text{times}}\hspace{-2mm},2}(1-z)(\sigma) ={\boldsymbol \zeta}_{{k}}^\ell(\sigma)
\quad (\sigma \in G_K).
\]
\end{nsTheorem}

\begin{nsRemark}
In \cite{S21}, the second named author showed that 
the functional equation (\ref{GaloisLi2FE}) has an application
to a reciprocity law of the triple mod-$\{2,3\}$ symbols
of rational primes via Ihara-Morishita theory
(cf. \cite{HM19}).
Theorem \ref{GaloisOiUeno} was shortly announced in 
a talk by the first named author at online Oberwolfach
meeting (\cite{N21}).
After the present paper was worked out, the second named author
obtained a generalization of Theorem \ref{GaloisOiUeno}
to higher multi-indices (\cite{S23a}), and showed
an $\ell$-adic version of Spence-Kummer's trilogarithm functional equation 
from which various formulas including Theorem \ref{l-adic Landen} can be derived 
by specializations 
(cf. \cite[Remark 4.3.]{S23b}).
\end{nsRemark}

The contents of this paper will be arranged as follows:
After a quick setup in \S 2 on the notations of standard paths on 
$\mathbf{P}^1-\{0,1,\infty\}$, in \S 3 we discuss complex
and $\ell$-adic Galois associators as formal power series in
two non-commuting variables, and define the multiple
polylogarithms as their coefficients of certain monomials.
We then review in the complex analytic context
that (\ref{LandenLi3}) and (\ref{OiUeno}) can be derived 
from algebraic relations (chain rules) of associators 
along simple compositions of paths. With this line in mind, 
we prove Theorems 1.1 and 1.2 in the $\ell$-adic
Galois case by tracing arguments in parallel ways to the complex case.
In \S 4, after shortly recalling polylogarithmic characters introduced in 
a series of collaborations by Wojtkowiak and the first named author, 
we present ${\mathbb Z}_\ell$-integrality test for $\ell$-adic Galois 
Landen's equation obtained in Theorem 1.1.
Section 5 turns to an alternative approach to functional equations of polylogarithms
based on a set of tools devised in \cite{NW12} to enhance 
Zagier's tensor criterion for functional equations into 
a concrete form.
Then we give alternative proofs of  (\ref{LandenLi3}) and Theorem 1.1 
with this method. 
In Appendix A,
%\ref{appendixA}, 
we exhibit lower degree terms of the 
complex and $\ell$-adic Galois associators
as a convenient reference from the text.
Appendix B
%\ref{appendixB} 
summarizes a set of computational tools 
from \cite{NW12} that converts a tensor criterion of functions into a 
polylogarithmic identity.

\bigskip
\noindent 
{\it Acknowledgement}: 
The authors would like to thank Hidekazu Furusho for valuable communications
on the subject
during the preparation of the present paper.
They also thank the anonymous referees for many pieces of advice 
which improved the presentation of this article.
This work was supported by JSPS KAKENHI Grant Numbers JP20H00115, JP20J11018.

\section{Set up}

Fix a prime number $\ell$.
Let $K$ be a subfield of the complex number field $\mathbb{C}$,
$\overline{K}$ the algebraic closure of $K$ in $\mathbb{C}$,
and $G_K:={\mathrm{Gal}}(\overline{K}/K)$
the absolute Galois group of $K$.
Let $U:=\mathbf{P}_{K}^1 - \{0,1,\infty\}$
 be the projective line minus three points over $K$,
$U_{\overline{K}}$ the base-change of $U$ via the inclusion $K\hookrightarrow \overline{K}$,
and
$U^{\rm an}={\mathbf P}^1({\mathbb C})- 
\{0,1,\infty\}$
the complex analytic space associated to the base-change of $U_{\overline{K}}$ 
via the inclusion $\overline{K} \hookrightarrow \mathbb{C}$.

In the following,
we shall write
$\overrightarrow{01}$
for the standard $K$-rational tangential base point on $U$.
Let $z$ be a $K$-rational point of $U$ or a $K$-rational tangential base point on $U$.
We consider $\overrightarrow{01}$,
$z$ also as points on $U_{\overline{K}}$ or $U^{\rm an}$ 
by inclusions $K \hookrightarrow \overline{K}$ and
$\overline{K} \hookrightarrow \mathbb{C}$.
({\it Note:}
We admit the particular case $K={\mathbb C}$ where $G_K=\{1\}$.
It is also possible to start with a specific complex point 
$z \in {\mathbb C}-\{0,1\}$ so that any field $K$ with 
${\mathbb Q}(z)\subset K\subset {\mathbb C}$ fits into our setup.)

Let 
$\pi_1^{\rm top}(U^{\rm an}; \overrightarrow{01}, z)$
be the set of homotopy classes of topological paths on $U^{\rm an}$ from $\overrightarrow{01}$ to $z$, and let   
%$\pi_1^\ellet(U_{\overline{K}}; \overrightarrow{01},{z})$
$\pi_1^{\ell\text{-\'et}}(U_{\overline{K}}; \overrightarrow{01},{z})$
be the pro-$\ell$-finite set of pro-$\ell$ \'etale paths on 
$U_{\overline{K}}$ from $\overrightarrow{01}$ to $z$.
Note that there is a canonical comparison map
\[
\pi_1^{\rm top}(U^{\rm an}; \overrightarrow{01}, z) \to \pi_1^{\ell\text{-\'et}}(U_{\overline{K}}; \overrightarrow{01},
{z})
\]
that allows us to consider 
topological paths on $U^{\rm an}$ as pro-$\ell$ \'etale paths on 
$U_{\overline{K}}$.

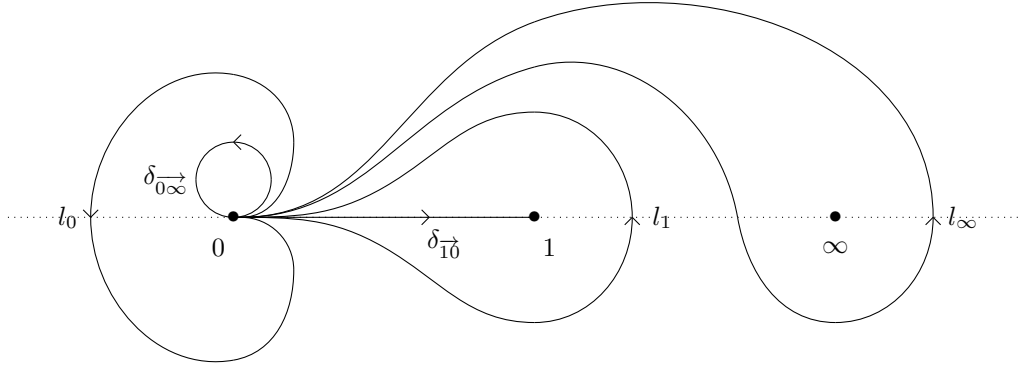
\begin{center}
%\framebox{
\begin{tikzpicture}[thick, scale=0.8]
\draw (2,0) -- (6,0);
\draw (4.6,0) -- (4.5,0.1);
\draw (4.6,0) -- (4.5,-0.1);

\draw (0.1,0) -- (0,0.1);
\draw (0.1,0) -- (0.2,0.1);

\draw (7.3,0) -- (7.4,-0.1);
\draw (7.3,0) -- (7.2,-0.1);

\draw (11.3,0) -- (11.4,-0.1);
\draw (11.3,0) -- (11.2,-0.1);

\draw (2,1) -- (2.1,1.1);
\draw (2,1) -- (2.1,0.9);

\draw (2.1,0) to [out=0,in=270] (2.8,1);
\draw (2.8,1) to [out=90,in=-10] (2,1.9);
\draw (2,1.9) to [out=170,in=90] (0.1,0);
\draw (0.1,0) to [out=270,in=190] (2,-1.9);
\draw (2,-1.9) to [out=10,in=270] (2.8,-0.7);

\draw (2.8,-0.7) to [out=90,in=360] (2.1,0);
\draw (2,0) to [out=0,in=210] (4,0.4);
\draw (4,0.4) to [out=30,in=180] (6,1.4);
\draw (6,1.4) to [out=0,in=90] (7.3,0);
\draw (7.3,0) to [out=270,in=0] (6,-1.4);
\draw (4,-0.4) to [out=330,in=180] (6,-1.4);
\draw (4,-0.4) to [out=150,in=0] (2,0);

\draw (2.1,0) to [out=0,in=195] (6,2);
\draw (6,2) to [out=15,in=100] (8.7,0);
\draw (8.7,0) to [out=280,in=180] (10,-1.4);
\draw (10,-1.4) to [out=0,in=270] (11.3,0);

\draw (2.1,0) to [out=0,in=200] (6,2.6);
\draw (6,2.6) to [out=20,in=90] (11.3,0);

\draw (2,0) to [out=0,in=270] (2.5,0.5);
\draw (2.5,0.5) to [out=90,in=0] (2,1);
\draw (2,1) to [out=180,in=90] (1.5,0.5);
\draw (1.5,0.5) to [out=270,in=180] (2,0);

\node at (-0.2,0) {$l_0$};
\node at (1.8,-0.4) {$0$};
\node at (4.8,-0.35) {$\delta_{\overrightarrow{10}}$};
\node at (1.1,0.5) {$\delta_{\overrightarrow{0\infty}}$};
\node at (6.2,-0.4) {$1$};
\node at (7.7,0) {$l_1$};
\node at (11.7,0) {$l_{\infty}$};
\node at (2,0) {${\bullet}$};
\node at (6,0) {${\bullet}$};
\node at (10,0) {${\bullet}$};
\node at (10,-0.4) {$\infty$};
\node at (6,-2.5) {The dashed line represents ${\mathbf P}^1({\mathbb R})- \{0,1,\infty\}$.};
\node at (6,-3) {The upper half-plane is above the dashed line.};

\draw[dotted](-1,0)--(12.5,0);
\end{tikzpicture}
%}
\end{center}

Let $l_0,
l_1,
l_\infty$
be the topological paths on $U^{\rm an}$
with base point
$\overrightarrow{01}$
circling counterclockwise around $0,
1,
\infty$,
respectively.
Then,
$\{l_0,
l_1\}$ is a free generating system of the topological fundamental group 
$\pi_1^{\rm top}(U^{\rm an}, \overrightarrow{01}):=\pi_1^{\rm top}(U^{\rm an}; \overrightarrow{01}, \overrightarrow{01})$ 
or the pro-$\ell$ \'etale fundamental group 
$\pi_1^{\ell\text{-\'et}}(U_{\overline{K}}, \overrightarrow{01}):=\pi_1^{\ell\text{-\'et}}(U_{\overline{K}}; \overrightarrow{01},
\overrightarrow{01})$.
Then, 
$\pi_1^{\rm top}(X^{\rm an}, \overrightarrow{01})$ is a free group of rank $2$ generated by 
$\{l_0,l_1\}$
and
$\pi_1^{\ell\text{-\'et}}(U_{\overline{K}}, \overrightarrow{01})$
is a free pro-$\ell$ group of rank $2$ topologically generated by $\{l_0,l_1\}$.
%\[
%\pi_1^\ellet(U_{\overline{K}}, \overrightarrow{01})=\overline{\left<l_0,
%l_1\right>}.
%\]

Fix a topological path
$\gamma_z \in \pi_1^{\rm top}(U^{\rm an}; \overrightarrow{01}, z)$
on $U^{\rm an}$ from $\overrightarrow{01}$ to $z$.
Moreover,
let
$\delta_{\overrightarrow{10}} \in \pi_1^{\rm top}(U^{\rm an}; {\overrightarrow{01}}, \overrightarrow{10})$
be the topological path on $U^{\rm an}$ 
from $\overrightarrow{01}$ to $\overrightarrow{10}$ along the real interval,
and let
$\delta_{\overrightarrow{0\infty}} \in \pi_1^{\rm top}(U^{\rm an}; {\overrightarrow{01}}, \overrightarrow{0\infty})$
be the topological path on the upper half-plane in $U^{\rm an}$ 
from $\overrightarrow{01}$ to $\overrightarrow{0\infty}$.

Let 
$\phi,
\psi \in {\rm Aut}(U^{\rm an})$
be
automorphisms of $U^{\rm an}$
defined by
\begin{equation}\label{S_3-sym}
\phi(t)=1-t,\quad
\psi(t)=\frac{t}{t-1}, 
\end{equation}
and introduce
specific paths from $\overrightarrow{01}$ to $1-z$ and to $\frac{z}{z-1}$ by
\begin{align}
\begin{dcases} \label{path-comp}
\gamma_{1-z}&:=\delta_{\overrightarrow{10}} \cdot \phi(\gamma_z) \in \pi_1^{\rm top}(U^{\rm an}; \overrightarrow{01}, 1-z), \\
\gamma_{\frac{z}{z-1}}&:=\delta_{\overrightarrow{0\infty}} \cdot \psi(\gamma_z) \in \pi_1^{\rm top}\left(U^{\rm an}; \overrightarrow{01}, \frac{z}{z-1}\right).
\end{dcases}
\end{align}
Here, paths are composed from left to right.

For any field $F$, we shall write 
$F\langle\!\langle X,Y\rangle\!\rangle$ for the ring of 
non-commutative power series 
in the non-commuting indeterminates $X,Y$
with coefficients in $F$.
Every element $f(X,Y)\in F\langle\!\langle X,Y\rangle\!\rangle$
can be expanded as a formal sum over the free
monoid ${\rm M}$ generated by $X,Y$.
We call an element $w\in {\rm M}$ 
a word in $X,Y$, and denote by $\mathtt{Coeff}_w(f(X,Y))\in F$
the coefficient of a word $w\in{\rm M}$
in $f(X,Y)\in F\langle\!\langle X,Y\rangle\!\rangle$:
\begin{equation}
\label{coeff-def}
f(X,Y)=\sum_{w\in {\rm M}}
\mathtt{Coeff}_w(f(X,Y))\cdot 
w.
\end{equation}
In particular, for the unit element $w=1$
of ${\rm M}$, $\mathtt{Coeff}_{1}(f(X,Y))$ denotes
the constant term of $f(X,Y)$.

\section{Associators and multiple polylogarithms}

Recall that the multiple polylogarithms
appear as coefficients of the non-commutative 
formal power series in two variables,
%with complex coefficients,
determined as the basic solution of the KZ equation (Knizhnik-Zamolodchikov equation) 
on $\mathbf{P}^1(\mathbb{C}) - \{0,1,\infty\}$
(cf. \cite{Dr90}).
More precisely, let
$G_0(X,Y)(z)\left(=G_0(X,Y)(\gamma_z)\right)$
be the fundamental solution
 of the formal KZ equation
\[
\dfrac{d}{d z}G(X,Y)(z)=
\left( \dfrac{X}{z}+\dfrac{Y}{z-1} \right)G(X,Y)(z)
\]
on $\mathbf{P}^1(\mathbb{C}) - \{0,1,\infty\}$,
which is an analytic function with values in
$\mathbb{C} \langle\!\langle X,
Y \rangle\!\rangle$
characterized by the asymptotic behavior
$G_0(X,Y)(\gamma_z) \approx \exp\left(X \log(z) \right)~(z \to 0)$
and analytically
continued to the universal cover of 
$\mathbf{P}^1({\mathbb C})-\{0,1,\infty\}$.
Here, $\log(z):=\int_{\delta_{\overrightarrow{01}}^{-1} \cdot \gamma_z}\frac{1}{t}dt$.

One can expand $G_0(X,Y)(\gamma_z)$ 
with the notation (\ref{coeff-def})
as:
\begin{equation}
G_0(X,Y)(\gamma_z)=1+\sum_{w \in {\rm M} \backslash \{1\}}
\mathtt{Coeff}_w(G_0(X,Y)(\gamma_z)) \cdot w.
\end{equation}
The multiple polylogarithm
$Li_\mathbf{k}(z)\left(=Li_\mathbf{k}(\gamma_z)\right)$
associated to a tuple 
$\mathbf{k}=(k_1\dots,k_d)$ $\in\mathbb{N}^d$ 
and a topological path $\gamma_z$ from $\overrightarrow{01}$ to $z$
is equal to the coefficient of $G_0(X,Y)(\gamma_z)$ at the `regular' word
$w(\mathbf{k}):=X^{k_d-1}Y\cdots X^{k_1-1}Y$
multiplied by $(-1)^d$ 
(where `regular' means that the word ends in the letter $Y$).
In summary, writing the length $d$ of the tuple $\mathbf{k}=(k_1\dots,k_d)$
as $\mathrm{dep}(\mathbf{k})$, we have
\begin{equation} \label{multiPoly}
Li_{\mathbf{k}}(\gamma_z)=(-1)^{\mathrm{dep}(\mathbf{k})}
\,
\mathtt{Coeff}_{w(\mathbf{k})}(G_0(X,Y)(\gamma_z)).
\end{equation}
To define the $\ell$-adic Galois multiple polylogarithms,
we make use of the $G_K$-action on the \'etale paths
instead of the fundamental KZ-solution.
Given a pro-$\ell$ \'etale path 
$\gamma_z \in \pi_1^{\ell\text{-\'et}}(U_{\overline{K}};\overrightarrow{01},z)$,
form a pro-$\ell$ \'etale loop
$
{\mathfrak f}^{\gamma}_{\sigma}:=\gamma \cdot \sigma(\gamma)^{-1} ~\in 
\pi_1^{\ell\text{-\'et}}(U_{\overline{K}}, \overrightarrow{01}),
$
and expand it via the Magnus embedding
$ \pi_1^{\ell\text{-\'et}}(U_{\overline{K}}, \overrightarrow{01}) \hookrightarrow
{\mathbb Q}_\ell\langle\!\langle X,Y\rangle\!\rangle
$ defined by $l_0\mapsto \exp(X)$,
$l_1\mapsto \exp(Y)$.

\begin{nsNotation}
We shall often identify  
${\mathfrak f}^{\gamma_z}_{\sigma} \in \pi_1^{\ell\text{-\'et}}(U_{\overline{K}},\overrightarrow{01})$
with the above image in ${\mathbb Q}_\ell\langle\!\langle X,Y\rangle\!\rangle$.
Then, following notation of (\ref{coeff-def}), let us 
write 
\begin{equation}
\label{Expansion_fsigma}
{\mathfrak f}^{\gamma_z}_{\sigma}(X,Y)=1+\sum_{w \in {\rm M} 
\backslash \{1\}}
\mathtt{Coeff}_{w}({\mathfrak f}^{\gamma_z}_{\sigma}(X,Y)) \cdot w
\qquad (\sigma\in G_K).
\end{equation}
\end{nsNotation}
\noindent
This is what we described in (\ref{Galois-associator-intro}).
In the parallel way to the above (\ref{multiPoly}),
for any tuple $\mathbf{k}$ of positive integers,
we define the $\ell$-adic Galois multiple polylogarithm
$Li_{\mathbf{k}}^\ell$ to be the function $G_K\to{\mathbb Q}_\ell$ 
determined by 
\begin{equation} \label{Gal_multiPoly}
Li_{\mathbf{k}}^\ell(\gamma_z)(\sigma)=(-1)^{\mathrm{dep}(\mathbf{k})}
\mathtt{Coeff}_{w(\mathbf{k})}({\mathfrak f}^{\gamma_z}_{\sigma}(X,Y)) 
\qquad (\sigma\in G_K).
\end{equation}
The $\ell$-adic zeta function
${\boldsymbol \zeta}_{{\mathbf k}}^\ell
:G_K \to {\mathbb Q}_\ell$
is a special case when $\gamma_z$ is the 
unit interval path 
$\delta_{\overrightarrow{10}}$
from $\overrightarrow{01}$ to $\overrightarrow{10}$
as mentioned in Introduction (\ref{l-adic_zeta}).

\begin{nsRemark}
It is worth noting that
the $\ell$-adic Galois associator
${\mathfrak f}^{\delta_{\overrightarrow{10}}}_{\sigma}(X,Y) 
\in \mathbb{Q}_{\ell} \langle\!\langle X, Y \rangle\!\rangle$
is the $\ell$-adic Galois analog of the
{Drinfeld} associator
\[
\Phi(X,Y):={\Bigl(G_0(Y,X)(\gamma_{1-z})\Bigr)}^{-1} \cdot G_0(X,Y)(\gamma_z) \in \mathbb{C} \langle\!\langle X,
Y \rangle\!\rangle.
\]
See also Appendix A
%\ref{appendixA} 
for some basic properties and explicit coefficients
of low degree terms of $\Phi(X,Y)$ and ${\mathfrak f}^{\delta_{\overrightarrow{10}}}_{\sigma}(X,Y)$.

\end{nsRemark}

\begin{nsLemma}[Key identities]
\label{KeyId} The notations being as above, the following identities hold.
\begin{enumerate} 
\item[(1)] $G_0(X,Y)(\gamma_z)=G_0(Y,X)(\gamma_{1-z}) \cdot \Phi(X,Y)$.
\item[(2)]  $G_0(X,Y)(\gamma_{\frac{z}{z-1}})=G_0(X,Z)(\gamma_z) \cdot {\rm{exp}}(\pi {\mathtt{i}} X),$
where $Z:=-X-Y$.
\item[(3)]  ${\mathfrak f}^{\gamma_{z}}_{\sigma}(X,Y)={\mathfrak f}^{\gamma_{1-z}}_{\sigma}(Y,
X) \cdot {\mathfrak f}^{\delta_{\overrightarrow{10}}}_{\sigma}(X,Y)$. 
\item[(4)] ${\mathfrak f}^{\gamma_{\frac{z}{z-1}}}_{\sigma}(X,Y)
   ={\mathfrak f}^{\gamma_{z}}_{\sigma}(X,Z) \cdot 
\exp(\frac{1-\chi(\sigma)}{2}X)$,
%   {\mathfrak f}^{\delta_{\overrightarrow{0\infty}}}_{\sigma}(X,Y),$ 
where $Z := {{\rm log}({{\rm exp}(-Y)}{{\rm exp}(-X)})}$.
\end{enumerate}
\end{nsLemma}

\begin{proof}
The identity (1) was remarked in \cite[A.24]{F14}. 
We shall prove (1) and (2) as consequences of the chain rule of iterated integrals along
composition of paths: $\Lambda(\alpha\beta)=\Lambda(\alpha) \Lambda(\beta)$ for
$\alpha:x{\leadsto} y$ and $\beta: y{\leadsto} z$, where
$\Lambda(\gamma)=1+\sum_{m=1}^\infty \int_\gamma \underbrace{\omega\cdots\omega}_{\scriptsize m}
\in {\mathbb C}\langle\!\langle X,Y\rangle\!\rangle$ is Chen's transport of the formal connection associated to 
$\omega=\frac{dt}{t}X+\frac{dt}{t-1}Y\in V_1\otimes \Omega^1$.
(Here $\Omega^1$ is the space of meromorphic 1-forms with log-singularities
on $(\mathbf{P}^1,\{0,1,\infty\})$ and $V_1$ is the dual of $\Omega^1$.
For extension to cases of tangential base points,
see \cite{De89}, \cite{W97}). 
Below, for simplicity, write $\delta:=\delta_{\overrightarrow{10}}$, ${\varepsilon}:=\delta_{\overrightarrow{0\infty}}$.
It is easy to see from (\ref{path-comp}) and $\phi(\delta)=\delta^{-1}$ that 
$\gamma_z=\delta\cdot \phi(\gamma_{1-z})$.
This identity and another  $\gamma_{\frac{z}{z-1}}={\varepsilon}\cdot\psi(\gamma_z)$ 
from  (\ref{path-comp}) 
imply
$\Lambda(\gamma_z)
=\Lambda(\delta) \Lambda(\phi(\gamma_{1-z}))
=\Lambda(\delta)\phi^{-1}_\ast(\Lambda(\gamma_{1-z}))$ and
$\Lambda(\gamma_{\frac{z}{1-z}})
=\Lambda({\varepsilon})\Lambda(\psi(\gamma_z))
=\Lambda({\varepsilon}) \psi^{-1}_\ast (\Lambda(\gamma_z))$
respectively (cf.~\cite[(4.5)]{NW12}). 
The assertions (1), (2) follow from them together with a normalization of convention: 
$G_0(X,Y)(\gamma_z)=\overline{\Lambda(\gamma_z)(X,Y)}$,
where $\overline{f(X,Y)}$ denotes the non-commutative 
formal power series in $X,Y$ obtained from $f(X,Y)$
by reversing the order of letters 
in each word (i.e., $\overline{x_1\cdots x_m\rule{0pt}{1.2ex}}=x_m\cdots x_1$ for 
$x_i\in\{X,Y\}$, $i=1,\dots,m$, $m\ge 1$).
Note that $\overline{\Lambda_1\Lambda_2}=\overline{\Lambda_2}\cdot\overline{\Lambda_1}$.

(3): From (\ref{path-comp}) again, we have $\gamma_z=\delta \cdot \phi(\gamma_{1-z})$
as above.
Noting that the automorphism $\phi$ is defined over $K$ and
that $\delta \phi(l_i) \delta^{-1} =l_{1-i}$ ($i=0,1$), we compute 
\begin{align*}
\mathfrak{f}_\sigma^{\gamma_z}
&=\gamma_z\cdot \sigma(\gamma_z)^{-1}   
=\delta\cdot \phi(\gamma_{1-z})\cdot \sigma(\phi(\gamma_{1-z})^{-1})
\cdot \sigma(\delta)^{-1} \\
&=\delta\cdot \phi(\mathfrak{f}_\sigma^{\gamma_{1-z}})\cdot
\delta^{-1} \cdot \mathfrak{f}_\sigma^{\overrightarrow{10}} 
= \mathfrak{f}_\sigma^{\gamma_{1-z}}(Y,X)\cdot  \mathfrak{f}_\sigma^{\overrightarrow{10}} (X,Y).
\end{align*}

(4): Recall that 
$l_\infty=\exp(Z)$ represents a loop in $\pi_1^{\ell\text{-\'et}}(U_{\overline{K}}, \overrightarrow{01})$
such that $l_0l_1l_\infty=1$. Noting that $\psi$ is defined over $K$
and that ${\varepsilon}\psi(l_0){\varepsilon}^{-1}=l_0$, ${\varepsilon}\psi(l_1){\varepsilon}^{-1}=l_\infty$, 
we compute
from 
$\gamma_{\frac{z}{z-1}}={\varepsilon} \cdot \psi(\gamma_z)$ (\ref{path-comp}):
\begin{align*}
\mathfrak{f}_\sigma^{\gamma_{\frac{z}{z-1}}}&=
{\varepsilon} \cdot \psi(\gamma_z)\cdot \sigma({\varepsilon} \cdot \psi(\gamma_z))^{-1} 
={\varepsilon}\cdot \psi(\mathfrak{f}_\sigma^{\gamma_z}) \cdot {\varepsilon}^{-1}
\cdot {\varepsilon} \sigma({\varepsilon})^{-1} \\
&=\mathfrak{f}_\sigma^{\gamma_z}(X,Z)
\cdot \mathfrak{f}_\sigma^{\varepsilon}(X,Y) 
=\mathfrak{f}_\sigma^{\gamma_z}(X,Z)
\cdot l_0^{-\frac{\chi(\sigma)-1}{2}}.
\end{align*}
In the last equality, we used a formula $\sigma({\varepsilon})=l_0^{\frac{\chi(\sigma)-1}{2}} {\varepsilon}$
$(\sigma\in G_K)$.
The proof of Lemma is completed.
\end{proof}

We summarize 
analogy between $\ell$-adic Galois
and complex associators as Table 1,
where the 3rd and 4th rows reflect the key identities 
of Lemma \ref{KeyId}.

\renewcommand{\arraystretch}{1.6}
\tabcolsep = 0.1cm
\begin{table}[htb]\label{table1}
\centering
  \caption{}
  \scalebox{0.9}{
  \begin{tabular}{|c||c|}  \hline
    $\ell$-adic Galois side & complex side  \\ \hline \hline
    ${\mathfrak f}^{\gamma_{z}}_{\sigma}(X,
Y) \in \mathbb{Q}_{\ell} \langle\!\langle X,Y \rangle\!\rangle$ & $G_0(X,Y)(\gamma_z) \in \mathbb{C} \langle\!\langle X,Y \rangle\!\rangle$ 
\\ \hline
     ${\mathfrak f}^{\delta_{\overrightarrow{10}}}_{\sigma}(X,Y) \in \mathbb{Q}_{\ell} \langle\!\langle X,
Y \rangle\!\rangle$ & $\Phi(X,Y) \in \mathbb{C} \langle\!\langle X,Y \rangle\!\rangle$  \\ \hline
  ${\mathfrak f}^{\gamma_{z}}_{\sigma}(X,Y)={\mathfrak f}^{\gamma_{1-z}}_{\sigma}(Y,
X) \cdot {\mathfrak f}^{\delta_{\overrightarrow{10}}}_{\sigma}(X,Y)$ 
 &  $G_0(X,Y)(\gamma_z)=G_0(Y,X)(\gamma_{1-z}) \cdot \Phi(X,Y)$
  \\ \hline
   ${\mathfrak f}^{\gamma_{\frac{z}{z-1}}}_{\sigma}(X,Y)
   ={\mathfrak f}^{\gamma_{z}}_{\sigma}(X,Z) \cdot 
   \exp(\frac{1-\chi(\sigma)}{2}X)$
%   {\mathfrak f}^{\delta_{\overrightarrow{0\infty}}}_{\sigma}(X,Y),$ 
   & $G_0(X,Y)(\gamma_{\frac{z}{z-1}})=G_0(X,Z)(\gamma_z) \cdot {\rm{exp}}(\pi {\mathtt{i}} X),$ \\
 $Z := {{\rm log}({{\rm exp}(-Y)}{{\rm exp}(-X)})}$ & $Z:=-Y-X$ \\ \hline
  ${Li}^\ell_{\mathbf{k}}(\gamma_z)(\sigma)$: $\ell$-adic {Galois} multiple polylog value & $Li_\mathbf{k}(\gamma_z)$: multiple polylog value \\ \hline 
   ${\boldsymbol \zeta}_{{\mathbf k}}^\ell(\sigma)$: $\ell$-adic {Galois} multiple zeta value & $\zeta(\mathbf{k})$: multiple zeta value  \\ \hline 
  \end{tabular}
  }
\end{table}

\noindent
{\bf Algebraic proof of (\ref{LandenLi3})-(\ref{OiUeno})}.
The following arguments are motivated from an enlightening remark given
in
Appendix of Furusho's lecture note \cite[A.24]{F14}.
By the explicit formula of Le-Murakami \cite{LM96} type due to 
Furusho 
\cite[Theorem 3.15]{F04}, the coefficient of $YX^{k-1}$ in $G_0(X,Y)(\gamma_z)$ is 
\begin{align*}
\mathtt{Coeff}_{YX^{k-1}}(G_0(X,Y)(\gamma_z))
&=
-\sum_{\substack{s+t=k-1 \\ s,t\ge 0}}
(-1)^s Li_{f'(Y\scalebox{0.6}[0.8]{\rotatebox[origin=c]{-90}{$\exists$}} X^s)}(z) \frac{\log^tz}{t!} \\
&=(-1)^k\sum_{t=0}^{k-1} (-1)^t Li_{k-t}(z)\frac{\log^t z}{t!},
\end{align*}
where $f'$ indicates the operation annihilating terms ending with
the letter $X$.
Applying this to the key identity 
$G_0(Y,X)(\gamma_{1-z})=G_0(X,Y)(\gamma_z) \cdot \Phi(Y,X)$ from 
Lemma \ref{KeyId} (1) (cf.~also \cite[A.24]{F14}), 
we see that
\begin{align*}
\mathtt{Coeff}_{YX^{k-1}}(G_0(Y,X)(\gamma_{1-z}))
&=
\mathtt{Coeff}_{XY^{k-1}}(G_0(X,Y)(\gamma_{1-z})) \\
&=(-1)^{k-1}
Li_{\hspace{-2mm}\underbrace{\scriptstyle 1,\ldots,1}_{k-2~\text{times}}\hspace{-2mm},2}~(1-z)
\end{align*}
is equal to 
\begin{align*}
&\mathtt{Coeff}_{YX^{k-1}}(G_0(X,Y)(\gamma_z))+\mathtt{Coeff}_{YX^{k-1}}(\Phi(Y,X)) \\
&=\mathtt{Coeff}_{YX^{k-1}}(G_0(X,Y)(\gamma_z))+\mathtt{Coeff}_{XY^{k-1}}(\Phi(X,Y)) \\
&=(-1)^k\sum_{t=0}^{k-1}(-1)^t Li_{k-t}(z)\frac{\log^tz}{t!}
+(-1)^{k-1}\zeta({\hspace{-0.5mm}\underbrace{1,\ldots,1}_{k-2~\text{times}},2})
%\hspace{-1mm}\underbrace{1,\ldots,1}_{k-2~\text{times}},2).
\end{align*}
Here we used %the symmetric relation
%$\Phi(X,Y)=\Phi(Y,X)^{-1}$ of the Drinfeld associator and 
a tautological identity 
$$
\mathtt{Coeff}_{w(X,Y)}(\Phi(X,Y))=\mathtt{Coeff}_{w(Y,X)}(\Phi(Y,X))
$$
and the fact that $\mathtt{Coeff}_{X^i}(\Phi(X,Y))=\mathtt{Coeff}_{Y^i}(\Phi(X,Y))=0$
for all $i\ge 1$.
This together with the well--known identity
$\zeta({\hspace{-0.5mm}\underbrace{1,\ldots,1}_{k-2~\text{times}},2})
%\zeta(\hspace{-1mm}\underbrace{1,\ldots,1}_{k-2~\text{times}},2)
=\zeta(k)$
 (duality formula)
derives the identity (\ref{OiUeno}). 

\begin{nsRemark}
The duality formula is known to be 
a consequence of the $2$-cycle relation and the shuffle product formula.
Cf.~\cite[Lemma 2.2]{F22},
\cite[p.12 (5)]{Sou13} or Appendix A
%\ref{appendixA} 
(\ref{duality-relation}).
\end{nsRemark}

Before going to prove (\ref{LandenLi3}), we compare the coefficients
of $YXY$ in the same identity 
$G_0(X,Y)(\gamma_z)=G_0(Y,X)(\gamma_{1-z})\cdot \Phi(X,Y)$ of
Lemma \ref{KeyId} (1).
For simplicity, we shall use the following abbreviated notation:
$$
c_w(\gamma_z):=\mathtt{Coeff}_w(G_0(X,Y)(\gamma_z))
$$ 
for any word $w\in{\rm M}$ and a path
$\gamma_z$ from $\overrightarrow{01}$ to a point $z$.
By simple calculation, we then obtain
\begin{align}
c_{YXY}(\gamma_z)&=-\zeta(2) c_{X}(\gamma_{1-z}) 
%c_{XY}(1)
 -2\zeta(3)%+c_{YXY}(1)
+c_{XYX}(\gamma_{1-z}) \\
&= -\zeta(2) c_{X}(\gamma_{1-z}) %c_{XY}(1)  
-2\zeta(3) %+c_{YXY}(1)
+\bigl(
c_{XY}(\gamma_{1-z})c_{X}(\gamma_{1-z})
-2c_{X^2Y}(\gamma_{1-z})\bigr)
\nonumber
\end{align}
where, in the former equality are used 
%倒置構文です identities は複数なので are を用いる
known identities 
$\mathtt{Coeff}_{XY}(\Phi(X,Y))=-\zeta(2)$, 
$\mathtt{Coeff}_{YXY}(\Phi(X,Y))=-2\zeta(3)$
(see Appendix A %\ref{appendixA} 
(\ref{Drinfeld-expansion}),(\ref{euler-rel})), 
and 
in the last equality is used the shuffle relation
%relation は単数なので isを用いる
according to $XY\scalebox{0.6}[0.8]{\rotatebox[origin=c]{-90}{$\exists$}} X=XYX+2X^2Y$
 (cf.~Appendix A %\ref{appendixA} 
 (\ref{example-shuffle})).
This leads to
\begin{equation}
\label{L_21} 
Li_{2,1}(z)=-\frac{\pi^2}{6}\log(1-z) -2\zeta(3) %+\zeta(2,1)
-Li_2(1-z)\log(1-z)+2Li_3(1-z).
\end{equation}
Now let us compare the coefficients of $X^2Y$ on both sides of the
key identity
$$
G_0(X,Y)(\gamma_{\frac{z}{z-1}})=G_0(X,Z)(\gamma_z) \cdot {\rm{exp}}(\pi {\mathtt{i}} X)
$$
from Lemma \ref{KeyId} (2). It follows easily that
\begin{equation} \label{c_XXY}
c_{XXY}(\gamma_{\frac{z}{z-1}})
=c_{XXY}(\gamma_z)
-c_{YYY}(\gamma_z)
+c_{XYY}(\gamma_z)
+c_{YXY}(\gamma_z),
\end{equation}
or equivalently,
\begin{equation}
\label{Li3zoverz-1}
-Li_3\left(\frac{z}{z-1}\right)=Li_3(z)+Li_{1,1,1}(z)+Li_{1,2}(z)+Li_{2,1}(z).
\end{equation}
We know from the case $k=3$ of (\ref{OiUeno}) with
interchange $z\leftrightarrow1-z$ that
\begin{equation}\label{L_12} 
Li_{1,2}(z)=\zeta(3)-\left(
Li_3(1-z)-Li_2(1-z)\log(1-z)-\frac12\log z \log^2(1-z)\right).
\end{equation}
Putting (\ref{L_21}) and (\ref{L_12}) into the last two terms of (\ref{Li3zoverz-1})
with noticing %some known low degree identities such as $\zeta(2,1)=-2\zeta(3)$, 
$Li_{1,1,1}(z)=-\frac16\log^3(1-z)$ (cf.\,Appendix A),
we obtain a proof of 
Landen's trilogarithm functional equation (\ref{LandenLi3}).
\qed 

\bigskip
\noindent
{\bf Proof of Theorem 1.2:}
In the $\ell$-adic Galois setting, the argument for the assertion goes 
in an almost parallel way to the above proof for (\ref{OiUeno}). 
In fact, the formula of Le-Murakami and Furusho type is generalized to
any group-like elements of ${\mathbb Q}_\ell\langle\!\langle X,Y\rangle\!\rangle$ in 
\cite{N23}, so that it holds that
\begin{equation}
\label{LeFurusho}
\mathtt{Coeff}_{YX^{k-1}}({\mathfrak f}^{\gamma_z}_{\sigma}(X,Y))=
%-\sum_{\substack{s+t=k-1 \\ s,t\ge 0}}
%(-1)^s Li_{f'(B\sha A^s)}(z) \frac{\log^tz}{t!}=
(-1)^k\sum_{t=0}^{k-1} Li_{k-t}^\ell(\gamma_z)\frac{\rho_z(\sigma)^t}{t!}.
\end{equation}
Comparing the coefficients of $YX^{k-1}$ in the key identity
$$
{\mathfrak f}^{\gamma_{z}}_{\sigma}(X,Y)
={\mathfrak f}^{\gamma_{1-z}}_{\sigma}(Y,X) \cdot {\mathfrak f}^{\delta_{\overrightarrow{10}}}_{\sigma}(X,Y)
$$ 
of Lemma \ref{KeyId} (3), we obtain
\begin{equation}
\label{l-adicOiUeno-pre}
\sum_{j=0}^{k-1}{Li}^\ell_{k-j}(\gamma_z)(\sigma)\frac{\rho_z(\sigma)^j}{j!} 
+{Li}^\ell_{\hspace{-2mm}\underbrace{\scriptstyle 1,\ldots,1}_{k-2~\text{times}}\hspace{-2mm},2}(\gamma_{1-z})(\sigma) =
{\boldsymbol \zeta}_{\hspace{-2mm}\underbrace{\scriptstyle 1,\ldots,1}_{k-2~\text{times}}\hspace{-2mm},2}^\ell(\sigma)
%{Li}^{\hspace{0.04cm}\ell,{\overrightarrow{10}}}_k(\sigma) 
\quad (\sigma \in G_K).
%\tag{$\#_{\ell.1}$},
\end{equation}
Note here that, 
in the special case $z=\overrightarrow{10}$ with $\gamma_z=\delta_{\overrightarrow{10}}$,
we should interpret that 
${Li}^\ell_{1,\dots,1,2}
(\gamma_{1-z})(\sigma) =0$ and that
$\rho_z(\sigma)^j=0,1$ according to whether $j>0$ or $j=0$,
from which we obtain the duality formula:
\begin{equation}
{\boldsymbol \zeta}_{\hspace{-2mm}\underbrace{\scriptstyle 1,\ldots,1}_{k-2~\text{times}}\hspace{-2mm},2}^\ell(\sigma)=
{\boldsymbol \zeta}_{k}^\ell(\sigma)\qquad(\sigma\in G_K).
\end{equation}
Putting this back to  (\ref{l-adicOiUeno-pre}) settles the proof of 
Theorem 1.2. \qed

\bigskip
\noindent
{\bf Proof of Theorem 1.1:}
We only have to examine the $\ell$-adic Galois versions
of the identities
(\ref{L_21}), (\ref{Li3zoverz-1}) and (\ref{L_12}) 
with replacing the role of $G_0(X,Y)(\gamma_\ast)$ 
by $\mathfrak{f}_\sigma^{\gamma_\ast}(X,Y)$.
It turns out that the two identities (\ref{L_21}), (\ref{L_12}) 
have exactly the parallel counterparts:
\begin{align}
\label{ladic_L_21} 
&Li^\ell_{2,1}(\gamma_z)(\sigma) \\
&={\boldsymbol \zeta}_2^\ell(\sigma) \rho_{1-z}(\sigma)
+{\boldsymbol \zeta}_{2,1}^\ell(\sigma)+Li^\ell_2(\gamma_{1-z})(\sigma)\rho_{1-z}(\sigma)
+2Li^\ell_3(\gamma_{1-z})(\sigma), \nonumber
\\
\label{ladic_L_12} 
&Li^\ell_{1,2}(\gamma_z)(\sigma) \\
&={\boldsymbol \zeta}^\ell_3(\sigma)
-\left(
Li^\ell_3(\gamma_{1-z})(\sigma)+Li^\ell_2(\gamma_{1-z})(\sigma)\rho_{1-z}(\sigma)
+\frac12\rho_z(\sigma) \rho_{1-z}(\sigma)^2 \right)
\nonumber 
\end{align}
with $\sigma\in G_K$.
There occurs a small difference for (\ref{Li3zoverz-1}) when
evaluating the key identity 
$$
{\mathfrak f}^{\gamma_{\frac{z}{z-1}}}_{\sigma}(X,Y)
={\mathfrak f}^{\gamma_{z}}_{\sigma}(X,Z) 
\cdot 
\exp\left(\frac{1-\chi(\sigma)}{2} X\right)
%{\mathfrak f}^{\delta_{\overrightarrow{0\infty}}}_{\sigma}(X,Y)
$$
from Lemma \ref{KeyId} (4)
with taking into accounts the Campbell-Hausdorff sum
\begin{align}
Z:=&{{\rm log}({{\rm exp}(-Y)}{{\rm exp}(-X)})}\\
=&-Y-X \quad
\underbrace{+\frac{1}{2}\left(YX-XY\right)-\frac{1}{12}XXY+\cdots}_{(\ast)}, \nonumber
\end{align}
where do exist nontrivial nonlinear terms $(\ast)$
in the $\ell$-adic Galois case.
To simplify our presentation, we shall use 
the following abbreviated notation:
$$
c_w^\ell(\gamma_z)(\sigma):=\mathtt{Coeff}_w({\mathfrak f}^{\gamma_z}_{\sigma}(X,Y))
$$ 
for any word $w\in{\rm M}$, a path
$\gamma_z$ from $\overrightarrow{01}$ to a point $z$
and $\sigma\in G_K$.
Then ${\mathfrak f}^{\gamma_z}_{\sigma}(X,Z)$ is calculated as follows:
\begin{align}
&{\mathfrak f}^{\gamma_z}_{\sigma}(X,Z)\\
=&1+c_{X}^{\ell}(\gamma_z)(\sigma) X+c_{Y}^{\ell}(\gamma_z)(\sigma)Z+c_{XY}^{\ell}(\gamma_z)(\sigma)XZ+c_{Y^2}^{\ell}(\gamma_z)(\sigma){Z}^2+\cdots 
\nonumber
\\
%=&1+\Coeff_{X}^{\ell}(\gamma_z)(\sigma) X+\Coeff_{Y}^{\ell}(\gamma_z)(\sigma)\left(-Y-X+\frac{1}{2}YX-\frac{1}{2}XY-\frac{1}{12}X^{2}Y+\cdots\right) \nonumber \\
%&+\Coeff_{XY}^{\ell}(\gamma_z)(\sigma)X\left(-Y-X+\frac{1}{2}YX-\frac{1}{2}XY+\cdots\right) \nonumber \\
%&+\Coeff_{Y^2}^{\ell}(\gamma_z)(\sigma)\left(-Y-X+\frac{1}{2}YX-\frac{1}{2}XY+\cdots\right)\left(-Y-X+\frac{1}{2}YX-\frac{1}{2}XY+\cdots\right) \nonumber \\
%&+\cdots \qquad (\sigma\in G_K).\nonumber \\
=&1+\left(c_{X}^{\ell}(\gamma_z)(\sigma)-c_{Y}^{\ell}(\gamma_z)(\sigma)\right) X \nonumber \\
&-c_{Y}^{\ell}(\gamma_z)(\sigma)Y \quad \underbrace{+c_{Y}^{\ell}(\gamma_z)(\sigma)\left(\frac{1}{2}YX-\frac{1}{2}XY-\frac{1}{12}XXY+\cdots\right)}_{\text{$\ell$-adic extra terms in $c_{Y}^{\ell}(\gamma_z)(\sigma)Z$}} \nonumber \\
&+\left(c_{Y^2}^{\ell}(\gamma_z)(\sigma)-c_{XY}^{\ell}(\gamma_z)(\sigma)\right)XX \quad \underbrace{+c_{XY}^{\ell}(\gamma_z)(\sigma)\left(\frac{1}{2}XYX-\frac{1}{2}XXY+\cdots\right)}_{\text{$\ell$-adic extra terms in $c_{XY}^{\ell}(\gamma_z)(\sigma)XZ$}} \nonumber \\
&+\left(c_{Y^2}^{\ell}(\gamma_z)(\sigma)-c_{XY}^{\ell}(\gamma_z)(\sigma)\right)XY \quad \underbrace{+c_{Y^2}^{\ell}(\gamma_z)(\sigma)\left(\frac{1}{2}XXY-\frac{1}{2}YYX+\cdots\right)}_{\text{$\ell$-adic extra terms in $\mathtt{Coeff}_{Y^2}^{\ell}(\gamma_z)(\sigma){Z}^2$}} \nonumber \\
&+c_{Y^2}^{\ell}(\gamma_z)(\sigma)YX
+c_{Y^2}^{\ell}(\gamma_z)(\sigma)YY+\cdots \qquad (\sigma\in G_K).  \nonumber
\end{align}
Summing up, we find
that the $\ell$-adic Galois analog to identity (\ref{c_XXY}) 
turns to get extra additional terms as:
\begin{align} \label{ladic_c_XXY} 
c^\ell_{XXY}(\gamma_{\frac{z}{z-1}})(\sigma)
&=-c^\ell_{XXY}(\gamma_z)(\sigma)
-c^\ell_{YYY}(\gamma_z)(\sigma)
+c^\ell_{XYY}(\gamma_z)(\sigma)+c^\ell_{YXY}(\gamma_z)(\sigma)
 \\
&\quad
-\left(
\frac{1}{2}c^\ell_{XY}(\gamma_z)(\sigma)
-\frac{1}{2}c^\ell_{YY}(\gamma_z)(\sigma)
+\frac{1}{12}c^\ell_{Y}(\gamma_z)(\sigma)
\right),
%\qquad (\sigma\in G_K),
\nonumber
\end{align}
from which follows that
\begin{align}
\label{ladic_Li3zoverz-1}  
Li^\ell_3(\gamma_{\frac{z}{z-1}})(\sigma)
&=-Li^\ell_3(\gamma_z)(\sigma)
-Li^\ell_{1,1,1}(\gamma_z)(\sigma)
-Li^\ell_{1,2}(\gamma_z)(\sigma)-Li^\ell_{2,1}(\gamma_z)(\sigma) \\
&\quad -\left(
\frac12 Li^\ell_2(\gamma_z)(\sigma)
+\frac14 \rho_{1-z}(\sigma)^2
+\frac{1}{12}\rho_{1-z}(\sigma)
\right)
%\qquad (\sigma\in G_K).
\nonumber
\end{align}
for $\sigma \in G_K$.
The asserted formula follows from (\ref{ladic_Li3zoverz-1}) after
$Li^\ell_{1,2}(\gamma_z)(\sigma)$, $Li^\ell_{2,1}(\gamma_z)(\sigma)$
in the RHS are replaced by the equations
 (\ref{ladic_L_12}), (\ref{ladic_L_21}) respectively
and from knowledge of the lower degree coefficients of 
$\mathfrak{f}_\sigma^{\gamma_\ast}(X,Y)$
illustrated 
in 
(\ref{lowerFsigma}), (\ref{example-duality}) and
(\ref{example-shuffle})
of Appendix A. 
%\ref{appendixA} 
%$\check$
%in $Y^2,Y^3$ (cf. Appendix) i.e. 
%$Li^\ell_{1,1}(\gamma_z)(\sigma)=\frac{1}{2}\rho_{1-z}(\sigma)^2$,
%$Li^\ell_{1,1,1}(\gamma_z)(\sigma)=\frac{1}{6}\rho_{1-z}(\sigma)^3$.
\qed
%%%%%%%%%%%%%%%%%%%%%%%%%%%%%%%%%
% References
%%%%%%%%%%%%%%%%%%%%%%%%%%%%%%%%%

\section{Polylogarithmic characters and ${\mathbb Z}_\ell$-integrality test}

There is a specific series of functions $\tilde{\chi}_{m}^{z}:G_K\to {\mathbb Z}_\ell$ 
(called the polylogarithmic characters) closely related to
the $\ell$-adic Galois polylogarithms $Li_k^\ell(z):G_K\to {\mathbb Q}_\ell$. 
Let $\gamma_z$ be an $\ell$-adic \'etale path from $\overrightarrow{01}$ to a $K$-rational (tangential) point $z$ on $\mathbb{P}^1-\{0,1,\infty\}$.

\begin{nsDefinition}[\cite{NW99}: $\ell$-adic Galois polylogarithmic character ]
\label{char}
For each $m \in {\mathbb N}$ and $\sigma \in G_K$,
we define $\tilde{\chi}_{m}^{\gamma_z}(\sigma)$ 
(often written shortly as $\tilde{\chi}_{m}^{z}(\sigma)$) 
by the (sequential) Kummer properties
\[
\zeta_{\ell^n}^{\tilde{\chi}_{m}^{z}(\sigma)}=
\sigma \left(\prod_{i=0}^{\ell^n -1}
(1-\zeta_{\ell^n}^{\chi(\sigma)
^{-1}i}z^{1/\ell^n})
^{\frac{i^{m-1}}{\ell^n}}\right) 
\bigg/ \prod_{i=0}^{\ell^n -1}
(1-\zeta_{\ell^n}
^{i+\rho_{z}(\sigma)}z^{1/\ell^n})
^{\frac{i^{m-1}}{\ell^n}}
\]
over $n\in {\mathbb N}$, where 
the roots $z^{1/n},
(1-z)^{1/n},
(1-\zeta^a_n z^{1/n})^{1/m}$
($n,m \in {\mathbb N},
a \in {\mathbb Z}$)
are chosen along the path $\gamma_z \in \pi_1^{\rm top}
(U^{\rm an}; \overrightarrow{01}, z)$,
$\rho_{z}(=\rho_{\gamma_z}): G_K \to {\mathbb Z}_{\ell}$
is the Kummer $1$-cocycle of the $\ell$-th power roots $\{z^{1/{\ell^n}}\}_{n}$ along $\gamma_z$,
and $\chi: G_K \to {\mathbb Z}_{\ell}^\times$ is the $\ell$-adic cyclotomic character.
We call the function
\[\tilde{\chi}_{m}
^{z} \left( =\tilde{\chi}_{m}
^{\gamma_z} \right) :
G_K \to {\mathbb Z}_{\ell}\]
{\it the ($\ell$-adic Galois) polylogarithmic character} associated to
$\gamma_z \in \pi_1^{\rm top}
(U^{\rm an}; \overrightarrow{01}, z)$.
When $z=\overrightarrow{10}$ and $\gamma_z=\delta_{\overrightarrow{10}}$,
it gives the Soul\'e character (cf. \cite[REMARK 2]{NW99}).
\end{nsDefinition}

We first begin with summarizing the relations between
the polylogarithmic characters and $\ell$-adic Galois polylogarithms:

\begin{nsProposition} \label{compareGaloisPolylog}
Let ${\mathfrak f}^{\gamma_z}_{\sigma}(X,Y)$ be the Magnus expansion
of the $\ell$-adic Galois associator 
${\mathfrak f}^{\gamma_z}_{\sigma}$
as in (\ref{Expansion_fsigma}).
Then, we have:
\begin{align*}
\mathtt{Coeff}_{YX^{m-1}}\left({\mathfrak f}^{\gamma_z}_{\sigma}(X,Y)\right)
&=-\frac{\tilde{\chi}_{m}^{z}(\sigma)}{(m-1)!}
%&
\ \left(=(-1)^m
\sum_{k=0}^{m-1}
{Li}^\ell_{m-k}(\gamma_z)(\sigma)
 \frac{\rho_z(\sigma)^k}{k!}
\right)
\tag{i}, \\
\mathtt{Coeff}_{X^{m-1}Y}\left({\mathfrak f}^{\gamma_z}_{\sigma}(X,Y)\right)
&
=-
{Li}^\ell_{m}(\gamma_z)(\sigma)
\left(=
(-1)^m \sum_{k=0}^{m-1} \frac{\rho_z(\sigma)^k}{k!}
\frac{\tilde{\chi}_{m-k}^{z}(\sigma)}{(m-1-k)!}
\right)
\tag{ii}
\end{align*}
for $\sigma\in G_K$.
\end{nsProposition}

\begin{proof}
%$\mathrm{dep}(YX^{k-1})=\mathrm{dep}(X^{k-1}Y)=1$
%and 
The first equality of (i) is proved in \cite[Proposition 8 (ii)]{NW20},
where the symbol $\mathtt{Li}_w$ in loc.cit. differs from our $Li_w$ by the
sign corresponding to the parity of the number of appearances of letter $Y$ in $w$.
The first equality of (ii) is just due to
our definition of $Li^\ell_{\mathbf{k}}$ in (\ref{Gal_multiPoly}).
The equality in the bracket of (i) is nothing but 
(\ref{LeFurusho}), i.e., is a consequence of 
a formula of Le-Murakami and Furusho type
(\cite{N23}).
% which compute
%the coefficient of any group-like element
%in $YX^{m-1}$ as an explicit linear combination of 
%the coefficients in $X^s Y$ ($0\le s \le m-1$).
The equality in the bracket of (ii) follows from (i) by
inductively reversing the sequence 
$\{\tilde{\chi}_{m}\}_m$ to
$\{ {Li}^\ell_{m}\}_m$.
\end{proof}

Often we prefer a functional equation of $\ell$-adic 
Galois polylogarithms converted
to 
a form of the
corresponding identity between
polylogarithmic characters
by Proposition \ref{compareGaloisPolylog}, because the latter enables us
to check the ${\mathbb Z}_\ell$-integrality of both sides of the equation.

For example, the functional equations (\ref{GaloisLi2FE}), 
(\ref{GaloisLi2Landen}) are respectively equivalent to 
\begin{align}
&\tilde{\chi}_{2}^{z}(\sigma)+\tilde{\chi}_{2}^{1-z}(\sigma)
+\rho_z(\sigma)\rho_{1-z}(\sigma)=\frac{1}{24}(\chi(\sigma)^2-1),
\label{GaloisTilchi2FE}
\\
&\tilde{\chi}_{2}^{z}(\sigma)+\tilde{\chi}_{2}^{z/(1-z)}(\sigma) 
=-\frac12 \rho_{1-z}(\sigma)(\rho_{1-z}(\sigma)-\chi(\sigma))
\label{GaloisTilchi2Landen}
\end{align}
for $\sigma\in G_K$. 
Noting that $\chi(\sigma)\equiv 1\pmod 2$ and 
$\chi(\sigma)^2\equiv 1\pmod {24}$, we easily see that each of the
RHSs have no denominator, i.e., $\in {\mathbb Z}_\ell$ for every prime $\ell$.
{}From this viewpoint, it is worth rewriting Landen's trilogarithm
functional equation (Theorem 1.1) in terms of polylogarithmic 
characters.
By simple computation, it results in:
\begin{align}
\label{Landen-Tri-characterform}
&\tilde{\chi}_{3}^{z}(\sigma) +\tilde{\chi}_{3}^{1-z}(\sigma)
+\tilde{\chi}_{3}^{z/(z-1)}(\sigma)
\\
&=
\tilde{\chi}_{3}^{\overrightarrow{10}}(\sigma)
+\chi(\sigma)\tilde{\chi}_{2}^{z}(\sigma)
+\rho_z(\sigma)\rho_{1-z}(\sigma)^2 
-\frac{\rho_{1-z}(\sigma)}{12}(\chi(\sigma)^2-1) 
\nonumber
\\
&\qquad -\frac{\rho_{1-z}(\sigma)}{6}
\Bigl(\chi(\sigma)-\rho_{1-z}(\sigma)\Bigr)
\Bigl(\chi(\sigma)-2\rho_{1-z}(\sigma)\Bigr).
\nonumber
\end{align}
It is not difficult to see that each term of
the above right--and side has
no denominator in ${\mathbb Q}_\ell$.

\section{Tensor criterion for Landen's equation for $Li_3$}

It would be worth giving alternative proofs of complex/$\ell$-adic Galois
Landen's trilogarithm functional equations (\ref{LandenLi3}) and Theorem 1.1 
with the method of \cite{NW12} not only for checking 
the validity of proofs given in \S 3 but also for providing a typical 
sample showing the utility of 
Zagier's tensor criterion for functional equations (cf.\,e.g.\,\cite{G13}).

Let $\mathcal{O}:={{\mathbb C}}[t,\frac{1}{t},\frac{1}{1-t}]$
be the coordinate ring of $U_{{\mathbb C}}=\mathbf{P}_{{\mathbb C}}^1 - \{0,1,\infty\}$ with
unit group ${\mathcal{O}}^{\times}$, and let
$f_1,
f_2,
f_3: U_{{\mathbb C}} \to U_{{\mathbb C}}$
be (auto)morphisms of $U_{{\mathbb C}}$ defined by
\[
f_1(t)=t,\quad
f_2(t)=1-t,\quad
f_3(t)=\frac{t}{t-1}.
\]
Considering $f_1,f_2,f_3:U\to \mathbf{G}_m$ as elements of $\mathcal{O}^{\times}$,
we specialize Zagier's tensor criterion for Landen's functional equation of $Li_3$'s
in the following proposition:

\begin{nsProposition}[Tensor criterion for Landen's functional equation for $Li_3$]\label{t-c}
Let $\overline{\mathcal{O}}^{\times} :={\mathcal{O}}^{\times} /{\mathbb C}^\times$ and denote
the image of $f_i$ $(i=1,2,3)$ in $\overline{\mathcal{O}}^{\times}$ by the same symbol. 
Then, in the tensor product 
$\overline{\mathcal{O}}^{\times} \otimes \bigl(\overline{\mathcal{O}}^{\times} \wedge 
\overline{\mathcal{O}}^{\times} \bigr)$
of multiplicative groups (where $\otimes$ and $\wedge$ are taken over $\mathbb{Z}$), we have
\[
f_1 \otimes (f_1 \wedge (f_1-1)) + f_2 \otimes (f_2 \wedge (f_2-1)) 
+ f_3 \otimes \left(f_3 \wedge \left(f_3-1\right)\right) \equiv 0.
\]
\end{nsProposition}
\begin{proof}
Set $a:=t$,
$b:=t-1$
$c:=-1$ and write the multiplication of $\mathcal{O}^{\times}$
in additive form. Then, we find $f_1=a$, $f_1-1=b$, 
$f_2=b+c$, $f_2-1=a+c$, $f_3=a-b$, $f_3-1=-b$ so that 
\begin{align*}
&f_1 \otimes (f_1 \wedge (f_1-1)) + f_2 \otimes (f_2 \wedge (f_2-1)) 
+ f_3 \otimes \left(f_3 \wedge \left(f_3-1\right)\right) \\
&= a \otimes (a \wedge b) +  (b+c) \otimes ((b+c) \wedge (a+c)) + (a-b) \otimes ((a-b) \wedge (-b))\\
&= b \otimes (c \wedge a) + b \otimes (b \wedge c) + c \otimes (a \wedge b) + c \otimes (a \wedge c) + c \otimes (b \wedge c).
\end{align*}
Obviously, this last side is annihilated in 
$\overline{\mathcal{O}}^{\times} \otimes \bigl(\overline{\mathcal{O}}^{\times} \wedge 
\overline{\mathcal{O}}^{\times} \bigr)$,
since $c\equiv 0$ in $\overline{\mathcal{O}}^{\times}.$
The assertion of proposition is proved.
\end{proof}

To compute the functional equations in concrete forms,  
we shall plug the above Proposition \ref{t-c} into
\cite[Remark 2.3 and Theorem 5.7: $\mathrm{(ii)}_{\mathbb C}\to \mathrm{(iii)}_{\mathbb C}$ and 
$\mathrm{(ii)}_\ell\to \mathrm{(iii)}_\ell$]{NW12}.
For a quick review on the method of \cite{NW12}, we also refer the reader 
to Appendix B. %\ref{appendixB}.
Fix a family of paths $\{\delta_1,\delta_2,\delta_3\}$ 
from $\overrightarrow{01}$
to 
$f_1(\overrightarrow{01})=\overrightarrow{01}$,
$f_2(\overrightarrow{01})=\overrightarrow{10}$,
$f_3(\overrightarrow{01})=\overrightarrow{0\infty}$,
with 
$\delta_1:=1(=\text{trivial path})$,
$\delta_2:=\delta_{\overrightarrow{10}}$,
$\delta_3:=\delta_{\overrightarrow{0\infty}}$
respectively.
Suppose we are given a topological path 
$\gamma_z:\overrightarrow{01}{\leadsto} z$ on 
$\mathbf{P}^1({\mathbb C})-\{0,1,\infty\}$.
Then, $\delta_i$ ($i=1,2,3$) provides a natural path 
$\delta_i \cdot f_i(\gamma_z):\overrightarrow{01}{\leadsto} f_i(z)$. 
Below we always consider the three points $f_1(z)=z$,
$f_2(z)=1-z$ and $f_3(z)=\frac{z}{z-1}$ to accompany
those natural tracking paths from the base point $\overrightarrow{01}$
$(i=1,2,3)$
in this way.

\subsection{Complex case}

With the notations being as above, 
\cite[Theorem 5.7 (iii$)_{\mathbb C}$]{NW12} asserts the existence of 
a functional equation of the form
\begin{equation} \label{complex-basic}
\sum_{i=1}^3
\mathcal{L}^{\varphi_3}_{{\mathbb C}}(f_i(z),f_i(\overrightarrow{01});f_i(\gamma_{z}))=0,
%\mathcal{L}^{\varphi_3}_{\C}(z,\overrightarrow{01};\gamma_{z})
%+\mathcal{L}^{\varphi_3}_{\C}(1-z,\overrightarrow{10};f_2(\gamma_{z}))
%+\mathcal{L}^{\varphi_3)}_{\C}\left(\frac{z}{z-1},\overrightarrow{0\infty};f_3(\gamma_{z})\right)=0,
\end{equation}
where each term can be calculated by a concrete algorithm 
\cite[Proposition 5.11]{NW12}.
Below let us exhibit the calculation by enhancing 
\cite[Examples 6.1-6.2]{NW12} to their 
``$Li_3$'' version, for which we start with the %$\log$-version of
graded Lie-versions of complex polylogarithms, written $\mathrm{li}_k(z,\gamma_z)$
for any path $\overrightarrow{01}{\leadsto} z$.
These can be converted to usual polylogarithms by 
\cite[Proposition 5.2]{NW12}; in particular, for $k=0,...,3$ we have: 
\begin{equation}
\label{liCtoLi}
\begin{dcases}
&\mathrm{li}_0(z,\gamma_z)=-\frac{1}{2\pi {\mathtt{i}}} \log(z), \\
&\mathrm{li}_1(z,\gamma_z)=-\frac{1}{2\pi {\mathtt{i}}} \log(1-z), \\
&\mathrm{li}_2(z,\gamma_z)=\frac{1}{4\pi^2} \left(Li_2(z)+\frac12 \log(z)\log(1-z)\right), \\
&\mathrm{li}_3(z,\gamma_z)= \frac{1}{(2\pi {\mathtt{i}})^3} \left(Li_3(z)-\frac12 \log(z) Li_2(z)
-\frac{1}{12} \log^2(z) \log(1-z)\right).
\end{dcases}
\end{equation}
Each term of (\ref{complex-basic}) relies only on the chain $f_i(\gamma_z)$
that does not start from $\overrightarrow{01}$ if $i\ne 1$, in which case
we need to interpret the chain $f_i(\gamma_z)$ as the difference 
``$\delta_i \cdot f_i(\gamma_z)$ minus $\delta_i$''.
At the level of graded Lie-version of polylogarithms, the difference 
can be evaluated by the polylog-BCH formula \cite[Proposition 5.9]{NW12}:
In our case, a crucial role is played by the polynomial
\begin{equation}
\label{PP3}
\mbox{\large $\mathsf{P}$}_3(\{a_j\}_{j=0}^3,\{b_j\}_{j=0}^3)
=a_3+b_3+\frac{1}{2}(a_0b_2-b_0a_2)
+\frac{1}{12}
(a_0^2b_1-a_0a_1b_0-a_0b_0b_1+a_1b_0^2)
\end{equation}
in 8 variables $a_j,b_j$ ($j=0,\dots,3$).
Using this and applying \cite[Proposition 5.11 (i)]{NW12}, we have
\begin{equation} \label{NW12Prop51}
\mathcal{L}^{\varphi_3}_{{\mathbb C}}(f_i(z),f_i(\overrightarrow{01});f_i(\gamma_{z}))
=\mbox{\large $\mathsf{P}$}_3\bigl(\{\mathrm{li}_j(f_i(z), \delta_i \cdot f_i(\gamma_z))\}_{j=0}^3,
\{-\mathrm{li}_j(f_i(\overrightarrow{01}),\delta_i)\}_{j=0}^3\bigr)
\end{equation}
for $i=1,2,3$.
Noting then that 
\begin{align*}
\left(-\mathrm{li}_j(\overrightarrow{01},\delta_1)\right)_{0\le j \le 3}&=(0,0,0,0), \\
\left(-\mathrm{li}_j(\overrightarrow{10},\delta_2)\right)_{0\le j \le 3} 
&=\left(0,0,-\mathrm{li}_2(\overrightarrow{10}),-\mathrm{li}_3(\overrightarrow{10})\right) \\
&=\left(0,0, -\frac{1}{4\pi^2}Li_2(1), -\frac{1}{(2\pi{\mathtt{i}})^3}Li_3(1)\right), \\
\left(-\mathrm{li}_j(\overrightarrow{0\infty},\delta_3)\right)_{0\le j\le 3}
&=\left(\frac{1}{2},0,0,0\right),
\end{align*}
we compute (\ref{NW12Prop51}) for $i=1,2,3$ as:
\begin{equation} 
\begin{dcases}
\mathcal{L}^{\varphi_3}_{{\mathbb C}}(z,\overrightarrow{01};\gamma_{z})
&=\mathrm{li}_3(z,\gamma_z), 
\\
\mathcal{L}^{\varphi_3}_{{\mathbb C}}(1-z,\overrightarrow{10};f_2(\gamma_{z}))
&=\mathrm{li}_3(1-z,\gamma_{1-z})-\mathrm{li}_3(\overrightarrow{10},\delta_{\overrightarrow{01}}) \\
&\quad+\tfrac12
\mathrm{li}_0(1-z,\gamma_{1-z})
(-\mathrm{li}_0(z,\gamma_z)), 
\\
\mathcal{L}^{\varphi_3}_{{\mathbb C}}\left(\tfrac{z}{z-1},\overrightarrow{0\infty};f_3(\gamma_{z})\right)
&=\mathrm{li}_3(\tfrac{z}{z-1},\gamma_{\frac{z}{z-1}})+\tfrac12\left(-\tfrac12 \mathrm{li}_2(\tfrac{z}{z-1},\gamma_{\frac{z}{z-1}})\right)\\
&\quad+\tfrac{1}{12}\left(\tfrac14 \mathrm{li}_1(\tfrac{z}{z-1},\gamma_{\frac{z}{z-1}})-\tfrac12\mathrm{li}_1(\tfrac{z}{z-1},\gamma_{\frac{z}{z-1}})
\mathrm{li}_0(\tfrac{z}{z-1},\gamma_{\frac{z}{z-1}})\right).
\end{dcases}
\end{equation}
Putting these together into (\ref{complex-basic}) and applying 
(\ref{liCtoLi}), we obtain Landen's functional equation (\ref{LandenLi3}).

\subsection{$\ell$-adic Galois case}
Let us apply \cite[Theorem 5.7 (iii$)_\ell$]{NW12}
in the parallel order to our above discussion in the complex case.
The $\ell$-adic version of the functional equation (\ref{complex-basic}) 
in loc.cit. relies on the choice of our free generator system 
$\vec{x}:=(l_0,l_1)$ of 
$\pi_1^{\ell\text{-\'et}}(U_{\overline{K}}, \overrightarrow{01})$
which plays an indispensable role to specify a splitting of the pro-unipotent
Lie algebra of $\pi_1^{\ell\text{-\'et}}$ into the weight gradation over ${\mathbb Q}_\ell$
(cf. \cite[\S 4.2]{NW12}).
Then, the functional equation turns out in the form 
\begin{equation} \label{ladic-basic}
\sum_{i=1}^3
\mathcal{L}^{\varphi_3(f_i),\vec{x}}_{\rm nv }(f_i(z),f_i(\overrightarrow{01});f_i(\gamma_{z}))(\sigma)=E(\sigma,\gamma_z)
\qquad (\sigma\in G_K)
\end{equation}
%in the terminologies of loc.\,cit., 
where $E(\sigma,\gamma_z)$ is called
the $\ell$-adic error term (\cite[\S 4.3]{NW12}).
The graded Lie-version of $\ell$-adic Galois polylogarithm
$\ell i_k(z,\gamma_z,\vec{x})$ (for $k\ge 1$)
is then defined as the coefficient 
of $\mathrm{ad}(X)^{k-1}(Y)=[X,[X,[\cdots [X,Y]..]$ 
in $\log({\mathfrak f}^{\gamma_z}_{\sigma}(X,Y)^{-1})$
as an element of Lie formal series $\mathrm{Lie}_{{\mathbb Q}_\ell}\langle\!\langle X,Y\rangle\!\rangle$.
Recall that the variables $X,Y$ are determined 
by $\vec{x}:=\{l_0,l_1\}$ by the
the Magnus embedding
$ \pi_1^{\ell\text{-\'et}}(U_{\overline{K}}, \overrightarrow{01}) \hookrightarrow
{\mathbb Q}_\ell\langle\!\langle X,Y\rangle\!\rangle
$ defined by $l_0\mapsto \exp(X)$,
$l_1\mapsto \exp(Y)$. 
For brevity below, let us often omit references to the loop system $\vec{x}=(l_0,l_1)$ 
and/or tracking paths $\delta_i \cdot f_i(\gamma_z):\overrightarrow{01}{\leadsto} f_i(z)$
in our notations as long as no confusions occur.
The list corresponding to (\ref{liCtoLi}) reads then:
\begin{equation} 
\label{li-GaltoLi}
\begin{dcases}
&\ell i_0(z,\gamma_z)(\sigma)=\rho_z(\sigma), \\
&\ell i_1(z,\gamma_z)(\sigma)=\rho_{1-z}(\sigma), \\
&\ell i_2(z,\gamma_z)(\sigma)=-\tilde{\chi}_2^z(\sigma)-\frac12\rho_z(\sigma)\rho_{1-z}(\sigma), \\
&\ell i_3(z,\gamma_z)(\sigma)=\frac12 \tilde{\chi}_3^z(\sigma)
+\frac12 \rho_z(\sigma) \tilde{\chi}_2^z(\sigma)
+ \frac{1}{12}\rho_z(\sigma)^2\rho_{1-z}(\sigma)
\end{dcases}
\end{equation}
with $\sigma\in G_K$.
Each term of the above (\ref{ladic-basic}) for $i=1,2,3$ can be expressed by 
the graded Lie-version of polylogarithms $\ell i_k$ ($k=0,\dots, 3$)
along ``$\delta_i \cdot f_i(\gamma_z)$ minus $\delta_i$''
by the polylog-BCH formula (\cite[Proposition 5.11 (ii)]{NW12})
in the following way:
\begin{align*}
&\mathcal{L}^{\varphi_3(f_i),\vec{x}}_{\rm nv }(f_i(z),f_i(\overrightarrow{01});f_i(\gamma_{z}))
%_\nv^{\varphi_n(f)_{\vec{x}}}(f(z),f(v);f(q))
%\tag*{(ii)} \\
%&
\\
&=
\mbox{\large $\mathsf{P}$}_3\bigl(
\{-\ell i_j(f_i(\overrightarrow{01}),\delta_i,\vec{x})\}_{j=0}^3,
\{\ell i_j(f_i(z),\delta_i \cdot f_i(\gamma_z),\vec{x})\}_{j=0}^3
\bigr).
\end{align*}
Noting that 
\begin{align*}
&\left(-\ell i_j(\overrightarrow{01},\delta_1)\right)_{0\le j \le 3}=(0,0,0,0), \\
&\left(-\ell i_j(\overrightarrow{10},\delta_2)\right)_{0\le j \le 3}
=\left(0,0,-\ell i_2(\overrightarrow{10}),-\ell i_3(\overrightarrow{10})\right)
=\left(0,0, \tilde{\chi}_2^{\overrightarrow{10}}(\sigma), -\frac{1}{2} \tilde{\chi}_3^{\overrightarrow{10}}(\sigma)\right), \\
&\left(-\ell i_j(\overrightarrow{0\infty},\delta_3)\right)_{0\le j\le 3}=\left(
\frac{1-\chi(\sigma)}{2},0,0,0\right),
\end{align*}
we obtain for $\sigma\in G_K$:
%[N-W2012]~Proposition 5.11, (5.3)
%\begin{eqnarray*}
\begin{align}
\label{LHSofGalBasic}
\begin{cases}
\mathcal{L}^{\varphi_3(f_1)}_{\rm nv}
&\hspace{-3mm}
(z,\overrightarrow{01};\gamma_{z})(\sigma) 
=\frac{1}{2} \tilde{\chi}_{3}^{z}(\sigma)+\frac{1}{2} \rho_{z}(\sigma)\tilde{\chi}_{2}^{z}(\sigma)
+\frac{1}{12} \rho_{z}(\sigma)^2 \rho_{1-z}(\sigma),
 \\
 & \\
\mathcal{L}^{\varphi_3(f_2)}_{\rm nv}
&\hspace{-3mm}(1-z,\overrightarrow{10};f_2(\gamma_{z}))(\sigma) \\
&\quad= \frac{1}{2} \tilde{\chi}_{3}^{1-z}(\sigma)
+\frac{1}{2} \rho_{1-z}(\sigma)\tilde{\chi}_{2}^{1-z}(\sigma)
+\frac{1}{12} \rho_{1-z}(\sigma)^2 \rho_{z}(\sigma)\\
& \qquad -\frac{1}{2}\tilde{\chi}_{3}^{\overrightarrow{10}}(\sigma)
-\frac{1}{2}\rho_{1-z}(\sigma)\tilde{\chi}_{2}^{\overrightarrow{10}}(\sigma), 
\\
& \\
\mathcal{L}^{\varphi_3(f_3)}_{\rm nv}
&\hspace{-3mm}
\left(\frac{z}{z-1},
\overrightarrow{0\infty};f_3(\gamma_{z})\right)(\sigma) \\
&\quad= 
\frac{1}{2} \tilde{\chi}_{3}^{\frac{z}{z-1}}(\sigma)+
\frac{1}{2} \rho_{\frac{z}{z-1}}(\sigma)\tilde{\chi}_{2}^{\frac{z}{z-1}}(\sigma)
+\frac{1}{12} \rho_{\frac{z}{z-1}}(\sigma)^2 \rho_{\frac{1}{1-z}}(\sigma)
\\
&\qquad +\frac{1}{2}\left( \frac{1-\chi(\sigma)}{2} \right) 
\left( -\tilde{\chi}_{2}^{\frac{z}{z-1}}(\sigma)-\dfrac{1}{2}\rho_{\frac{z}{z-1}}(\sigma)
\rho_{\frac{1}{1-z}}(\sigma) \right)
\\
&\qquad  +\frac{1}{12}\left( \frac{1-\chi(\sigma)}{2} \right)^2 \rho_{\frac{1}{1-z}}(\sigma)
-\frac{1}{12}\left( \frac{1-\chi(\sigma)}{2} \right)\rho_{\frac{z}{z-1}}(\sigma)\rho_{\frac{1}{1-z}}(\sigma).
%\end{eqnarray*}
\end{cases}
\end{align}
Combining the identities in (\ref{LHSofGalBasic}) enables us to 
rewrite the LHS of (\ref{ladic-basic}) in terms of $\ell$-adic Galois
polylogarithmic characters. 
It remains to compute the error term $E(\sigma,\gamma_z) $
in the right--and side of (\ref{ladic-basic}).

\begin{nsLemma}
\label{error-term-Landen}
Notations being as above, we have
\[
E(\sigma,\gamma_z)=-\dfrac{1}{12}\rho_{1-z}(\sigma)+\dfrac{1}{2} {\tilde{\chi}_{2}}
^{z}(\sigma)+\dfrac{1}{4}\rho_{z}(\sigma)\rho_{1-z}(\sigma).
\]
\end{nsLemma}
\begin{proof}
We shall apply the formula \cite[Corollary 5.8]{NW12} to compute the error term.
Let $[{\rm log}({\mathfrak f}^{\gamma_z}_{\sigma})^{-1}]_{<3}$ be 
the part of degree $<3$ cut out from the Lie formal series 
${\rm log}({\mathfrak f}^{\gamma_z}_{\sigma})^{-1}
\in \mathrm{Lie}_{{\mathbb Q}_\ell}\langle\!\langle X,Y\rangle\!\rangle$ with respect
obtained by the Magnus embedding 
$l_0\to e^X$, $l_1\to e^Y$ with respect 
to the fixed free generator system $\vec{x}=(l_0,l_1)$ of 
$\pi_1^{\ell\text{-\'et}}(U_{\overline{K}}, \overrightarrow{01})$.
We also write $\varphi_3: \mathrm{Lie}_{{\mathbb Q}_\ell}\langle\!\langle X,Y\rangle\!\rangle\to {\mathbb Q}_\ell$
for the ${\mathbb Q}_\ell$-linear form that picks up the coefficient of 
$[X,[X,Y]]$ (that is uniquely determined) for any Lie series of 
$\mathrm{Lie}_{{\mathbb Q}_\ell}\langle\!\langle X,Y\rangle\!\rangle $.
Introduce the variable $Z$ so that $e^Xe^Ye^Z=1$ in ${\mathbb Q}_\ell\langle\!\langle X,Y\rangle\!\rangle$.
By the Campbell-Baker-Hausdorff formula, we have
\begin{equation} \label{Z_BCH}
Z=-X-Y-\frac{1}{2}[X,Y]
-\dfrac{1}{12}[X,[X,Y]]+\cdots.
\end{equation}
According to \cite[Corollary 5.8]{NW12}, it follows then that
\begin{eqnarray*}
E(\sigma,\gamma_z) &=&
\sum_{i=1}^3 \varphi_3\left(\delta_i \cdot f_i
\left([{\rm log}({\mathfrak f}^{\gamma_z}_{\sigma})^{-1}]_{<3}\right) \cdot {\delta_i}^{-1}\right)
\\
&=& \sum_{i=1}^3 \varphi_3\Bigl(
\delta_i \cdot f_i
\left(\rho_{z}(\sigma)X+\rho_{1-z}(\sigma)Y+
\ell i_{2}(z,\gamma_z)(\sigma)[X,Y]\right) \cdot {\delta_i}^{-1}
\Bigr)\\
&=& \varphi_3\Bigl(\rho_{z}(\sigma)X+\rho_{1-z}(\sigma)Y
+\ell i_{2}(z,\gamma_z)(\sigma)[X,Y]\Bigr) \\
&& +\varphi_3\Bigl(\rho_{z}(\sigma)Y+\rho_{1-z}(\sigma)X
+\ell i_{2}(z,\gamma_z)(\sigma)[Y,X]\Bigr) \\
&& +\varphi_3\Bigl(\rho_{z}(\sigma)X+\rho_{1-z}(\sigma)Z
+\ell i_{2}(z,\gamma_z)(\sigma)[X,Z]\Bigr).
% \\
%&=& \varphi_3\left(\rho_{z}(\sigma){\rm log}(\overline{l_0})+\rho_{1-z}(\sigma){\rm log}%(\overline{l_\infty})+\ell i_{2}(z,\gamma_z)(\sigma)[{\rm log}(\overline{l_0}),{\rm log}(\overline{l_\infty})]\right) \\
\end{eqnarray*}
Here in the last equality, we applied the following table where 
$\delta_i \cdot f_i(\#) \cdot \delta_i^{-1}$
for $i \in \{1,2,3\}$, $\# \in \{X,Y\}$ are summarized.
\renewcommand{\arraystretch}{1.8}
\tabcolsep = 0.5cm
\begin{table}[hbtp]
  %\caption{}
  %\label{}
  \centering
  \begin{tabular}{|c||c|c|c|}
    \hline
    $i$ & $1$ & $2$ & $3$ \\
    \hline \hline
    $\delta_i \cdot f_i(X) \cdot \delta_i^{-1}$ & $X$ & $Y$ & $X$ \\
    \hline 
    $\delta_i \cdot f_i(Y) \cdot \delta_i^{-1}$  & $Y$ & $X$  & $Z$ \\ 
    \hline  
  \end{tabular}
\end{table}

\noindent
Since $\varphi_3$ annihilates those terms $X,Y,[X,Y], [Y,X]$, we continue the above
computation after (\ref{Z_BCH}) as:
\begin{eqnarray*}
E(\sigma,\gamma_z) &=& \varphi_3\left( -\dfrac{1}{12}\rho_{1-z}(\sigma)
[X,[X,Y]]-\dfrac{1}{2} \ell i_{2}(z,\gamma_z)(\sigma) [X,[X,Y]] \right)\\
&=& -\dfrac{1}{12}\rho_{1-z}(\sigma) -\dfrac{1}{2} \ell i_{2}(z,\gamma_z)(\sigma)\\
&=& -\dfrac{1}{12}\rho_{1-z}(\sigma) -\dfrac{1}{2}\left(-{\tilde{\chi}_{2}}
^{z}(\sigma)-\dfrac{1}{2}\rho_{z}(\sigma)\rho_{1-z}(\sigma)\right)\\
&=& -\dfrac{1}{12}\rho_{1-z}(\sigma)+\dfrac{1}{2} {\tilde{\chi}_{2}}
^{z}(\sigma)+\dfrac{1}{4}\rho_{z}(\sigma)\rho_{1-z}(\sigma).
\end{eqnarray*}
This concludes the assertion of the lemma.
\end{proof}

\begin{proof}[Alternative proof of Theorem 1.1]
As discussed in \S 4, 
the $\ell$-adic Galois Landen's trilogarithm functional equation in Theorem 1.1
is equivalent to the identity (\ref{Landen-Tri-characterform})
between polylogarithmic characters. 
The latter follows from (\ref{ladic-basic}) with replacements of the terms of both sides
by (\ref{LHSofGalBasic}) and Lemma \ref{error-term-Landen} by simple computations.
\end{proof}

%\appendix\section*{Appendix}
\section*{Appendix A: Low degree terms of associators}
\label{appendixA}
Presentation of lower degree terms of $G_{0}(X,Y)(z)$ and 
${\mathfrak f}^{\gamma_z}_{\sigma}(X,Y)$ are often useful as references.
The former one presented below reconfirms Furusho's preceding computations
found in \cite[3.25]{F04}-\cite[A.16]{F14} (where the sign of $\log (z)Li_2(z)$ 
had an unfortunate misprint in the coefficient of $XYX$).

\begin{align}
%$\displaystyle 
G_{0}&(X,Y)(z)=1+{\rm log}(z)X+
{\rm log}(1-z)Y+\dfrac{{\rm log}^{2}(z)}{2}{X}^2-{Li}_{2}(z)XY 
\\
&+\Bigl({Li}_{2}(z)+{\rm log}(z){\rm log}(1-z) \Bigr)YX+\dfrac{{{\rm log}^{2}(1-z)}}{2}Y^{2}
+\dfrac{{\rm log}^{3}(z)}{6}{X}^3-{Li}_3(z)X^{2}Y
\nonumber \\
&+\Bigl(2{Li}_{3}(z)-{\rm log}(z){Li}_{2}(z)\Bigr)XYX+{Li}_{1,2}(z)XY^{2}
\nonumber \\
&-\left( {Li}_{3}(z)-{\rm log}(z){Li}_{2}(z)-\dfrac{{\rm log}^{2}(z){\rm log}(1-z)}{2} \right)
YX^{2}+{Li}_{2,1}(z)YXY 
\nonumber \\
&-\left( {Li}_{1,2}(z)+{Li}_{2,1}(z)
-\dfrac{{\rm log}^{}(z){\rm log}^{2}(1-z)}{2} \right)Y^{2}X+\dfrac{{\rm log}^{3}(1-z)}{6}{Y}^3
\nonumber \\
&+\cdots \text{(higher degree terms)}.\nonumber
\end{align}
This is a group-like element of ${\mathbb C}\langle\!\langle X,Y\rangle\!\rangle$ whose coefficients
satisfy what are called the shuffle relations (\cite{Ree58}).
The regular coefficients (viz. those coefficients of monomials ending with the letter $Y$)
are given by iterated integrals of a sequence of $dz/z$, $dz/(1-z)$.
This immediately shows $G_0(0,Y)(z)=\sum_{k=0}\frac{\log^k(1-z)}{k!} Y^k$ and
say, $Li_{1,1,1}(z)=-\frac16\log^3(1-z)$.
Furusho gave an explicit formula that expresses arbitrary coefficients of
$G_0(X,Y)$ in terms only of the regular coefficients (\cite[Theorem 3.15]{F04}).
The specialization $z\to\overrightarrow{10}$ with $\gamma_z=\delta_{\overrightarrow{10}}$
(cf. \cite{W97} p.239 for a naive account)
interprets $\log z\to 0, \log(1-z)\to 0$ so as to produce the Drinfeld's associator:
\begin{align}
\label{Drinfeld-expansion}
&\Phi(X,Y)\quad \left(=G_0(X,Y)(\overrightarrow{10})\right) 
\\
=&~1-\zeta(2)XY+\zeta(2)YX 
-\zeta(3)X^2Y + 2\zeta(3)XYX\nonumber
 \\
&+\zeta(1,2)XY^2-\zeta(3)YX^2-2\zeta(1, 2)YXY
+\zeta(1,2)Y^2X \nonumber
\\
&+\cdots \text{(higher degree terms)}\nonumber 
\nonumber
\end{align}
which plays a primary role to define the Grothendieck-Teichm\"uller group
(\cite{Dr90}, \cite{Ih90}).

The expansion in ${\mathbb Q}_\ell\langle\!\langle X,Y\rangle\!\rangle$
of the $\ell$-adic Galois associator ${\mathfrak f}^{\gamma_z}_{\sigma}\in
\pi_1^{\ell\text{-\'et}}(U_{\overline{K}}, \overrightarrow{01})$
via the Magnus embedding $l_0\mapsto e^X$, $l_1\mapsto e^Y$ over ${\mathbb Q}_\ell$ reads as follows:

{\allowdisplaybreaks
{\small
\begin{align}
\label{lowerFsigma}
%$\displaystyle 
&{\mathfrak f}^{\gamma_z}_{\sigma}(X,Y)
=1-\rho_{\gamma_z}(\sigma)X-\rho_{\gamma_{1-z}}(\sigma)Y
+\dfrac{{\rho_{\gamma_z}(\sigma)}^2}{2}{X}^2-{Li}^\ell_{2}(\gamma_z)(\sigma)XY 
 \\
&+\Bigl( {Li}^\ell_{2}(\gamma_z)(\sigma)+\rho_{\gamma_z}(\sigma)\rho_{\gamma_{1-z}}(\sigma) \Bigr)YX+\dfrac{{\rho_{\gamma_{1-z}}(\sigma)}^2}{2}{Y}^2 -\dfrac{{\rho_{\gamma_z}(\sigma)}^3}{6}{X}^3-{Li}^\ell_{3}(\gamma_z)(\sigma)X^{2}Y 
 \nonumber \\
&+\Bigl( 2{Li}^\ell_{3}(\gamma_z)(\sigma)+\rho_{\gamma_z}(\sigma){Li}^\ell_{2}(\gamma_z)(\sigma) \Bigr)XYX+{Li}^\ell_{1,2}(\gamma_z)(\sigma)XY^{2} 
\nonumber \\
&-\left( {Li}^\ell_{3}(\gamma_z)(\sigma)+\rho_{\gamma_z}(\sigma){Li}^\ell_{2}(\gamma_z)(\sigma)+\dfrac{{\rho_{\gamma_z}(\sigma)}^{2}\rho_{\gamma_{1-z}}(\sigma)}{2} \right)YX^{2}+{Li}^\ell_{2,1}(\gamma_z)(\sigma)YXY 
\nonumber \\
&-\left( {Li}^\ell_{1,2}(\gamma_z)(\sigma)+{Li}^\ell_{2,1}(\gamma_z)(\sigma)+\dfrac{{\rho_{\gamma_z}(\sigma)}{\rho_{\gamma_{1-z}}(\sigma)}^{2}}{2} \right)Y^{2}X-\dfrac{{\rho_{\gamma_{1-z}}(\sigma)}^3}{6}{Y}^3 
\nonumber \\
&+\cdots \text{(higher degree terms)} 
\qquad (\sigma\in G_K).\nonumber
\end{align}
}
}

\noindent
The coefficients of $X,Y,YX^{k}$ $(k=1,2,...)$ were calculated in terms of
polylogarithmic characters explicitly in \cite{NW99}.
A formula of Le-Murakami, Furusho type for arbitrary group-like power series
was shown in \cite{N23}. As illustrated in Proposition \ref{compareGaloisPolylog},
the family of polylogarithmic characters and that of $\ell$-adic Galois polylogarithms are converted to each other.
The terms appearing in the above (\ref{lowerFsigma}) can be derived 
from them.

The $\ell$-adic Galois associator $\mathfrak{f}_\sigma^{\overrightarrow{10}}(X,Y)$ specialized 
at $z=\overrightarrow{10}$ with $\gamma_z=\delta_{\overrightarrow{10}}$ plays an important role
to define the pro-$\ell$ version of the Grothendieck-Teichm\"uller group 
similarly to the Drinfeld associator $\Phi(X,Y)$.

Let $f(X,Y)\in F\langle\!\langle X,Y\rangle\!\rangle$ 
be one of the power series either of  
$\Phi(X,Y)\in{\mathbb C}\langle\!\langle X.Y\rangle\!\rangle$ or
$\mathfrak{f}_\sigma^{\overrightarrow{10}}(X,Y)\in {\mathbb Q}_\ell\langle\!\langle X,Y\rangle\!\rangle$ 
for some $\sigma\in G_{\mathbb Q}$ with coefficients in $F={\mathbb C},{\mathbb Q}_\ell$
respectively. The following properties are well known to be held by $f(X,Y)$:

\begin{minipage}{0.9\textwidth}
\begin{enumerate}
\item[(i)] $f(X,Y)$ is a group-like element, i.e., the coefficients 
satisfy the shuffle relations (\cite{Ree58});
\item[(ii)] $f(X,0)=f(0,Y)=f(X,X)=1$ and 
$f(X,Y)\equiv 1$ mod $(X,Y)^2$;
\item[(iii)] $f(X,Y)f(Y,X)=1$ (2-cyclic relation).
\end{enumerate}
\end{minipage}

\smallskip
\noindent
The condition (iii) combined with (i) is also known to be equivalent to the following
\\
{\bf Duality-relation:}
If a power series 
$\displaystyle f(X,Y)=1+\sum_{w\in\mathrm{M}}c_w w\in F\langle\!\langle X,Y\rangle\!\rangle$ 
 $(c_w\in F)$ satisfies the above conditions (i) and (iii), then
\begin{equation}
\label{duality-relation}
c_w=(-1)^{|w|} c_{\overline{w'}}.
\end{equation}
Here, for a word $w=x_1\cdots x_m$ ($x_i=X,Y$), we write
$|w|=m$, and designate $\overline{w'}:=\overline{x_m'\cdots x_1'}$ to mean the word 
obtained by applying the substitutions ${X}'=Y$, ${Y}'=X$ 
after reversing the order of letters in $w$.
\\
See, e.g., discussions around \cite[p.12 (5)]{Sou13} for an algebraic proof of the duality relation
under (i), (iii).
The duality relation is already useful in degree 3, for example, to obtain
\begin{equation} 
\label{example-duality}
c_{XXY}=-c_{XYY}, \quad
c_{YXY}=-c_{XYX}, \quad
c_{YXX}=-c_{YYX}.
\end{equation} 
The first equation implies Euler's celebrated relation 
\begin{equation}
\label{euler-rel}
\zeta(3)=\zeta(1,2)
\end{equation}
when $F={\mathbb C}$ and its $\ell$-adic Galois analog 
${\boldsymbol\zeta}^\ell(3)(\sigma)={\boldsymbol\zeta}^\ell(1,2)(\sigma)$
when $F={\mathbb Q}_\ell$ (recall (\ref{l-adic_zeta}) for the latter notation).
One also observes that the shuffle relations corresponding to $XY\scalebox{0.6}[0.8]{\rotatebox[origin=c]{-90}{$\exists$}} Y=YXY+2XYY$,
$Y\scalebox{0.6}[0.8]{\rotatebox[origin=c]{-90}{$\exists$}} XX=YXX+XYX+XXY$ imply 
\begin{equation}
\label{example-shuffle}
c_{YXY}+2c_{XYY}=0=c_{YXX}+c_{XYX}+c_{XXY}.
\end{equation}
These equalities in (\ref{example-duality}), (\ref{example-shuffle})
enable us to find the above description of the degree 3 part of 
(\ref{Drinfeld-expansion}) and its obvious $\ell$-adic Galois analog.

\section*{Appendix B: A quick review of \cite{NW12} with remarks }
\label{appendixB}

Let $U=\mathbf{P}^1-\{0,1,\infty\}$ and let
$\pi=\pi_1(U_{\mathbb C},\overrightarrow{01})$ be the discrete fundamental group.
Let $m\ge 2$ be an integer, and 
suppose we are given a ${\mathbb Z}$-homomorphism 
$\varphi:\mathrm{gr}^m(\pi)\to {\mathbb Z}$ from the $m$-th 
graded quotient of the lower central series of $\pi$.
Let $H$ be the abelianization $\pi/[\pi,\pi]$ of $\pi$,
and define the standard projection
\begin{equation}
\mathrm{st}:
H^{\otimes m-2}\otimes \wedge^2 H
\twoheadrightarrow
\mathrm{gr}^m(\pi)
\end{equation}
by %x_1\otimes \cdots \otimes x_m \mapsto 
$
x_1\otimes \cdots \otimes x_{m-2}
\otimes (x_{m-1}\wedge x_m)
\mapsto 
[x_1,[x_2,\dots,[x_{m-1},x_m]..]]
$.
Let $\mathcal{O}={\mathbb C}[t,\frac{1}{t},\frac{1}{t-1}]$ be the
affine ring defining $U_C$.
In \cite[Corollary 3.7]{NW12}, it is shown that there
is an element 
\begin{equation}
\widehat{\kappa_{\otimes m}}(\varphi)\in
{{\mathcal{O}}^{\times}}^{\otimes n-2} 
\otimes 
\bigl({\mathcal{O}}^{\times} \wedge 
{\mathcal{O}}^{\times} 
\bigr)
\end{equation}
whose multi-Kummer dual 
$\kappa^{\otimes m}(\widehat{\kappa_{\otimes m}}(\varphi))
: H^{\otimes m-2}\otimes \wedge^2(H)\to {\mathbb Z}$
coincides with the composite $\varphi\circ \mathrm{st}$.
Note here that $\wedge^2H$ is understood as the
quotient wedge tensor (the maximal quotient of $H\otimes H$
satisfying $x\wedge y+y\wedge x=0$), while 
${\mathcal{O}}^{\times} \wedge 
{\mathcal{O}}^{\times} $ is understood as 
the submodule wedge tensor (the submodule generated by
the elements of the form $a\otimes b-b\otimes a$ in 
$({\mathcal{O}}^{\times} )^{\otimes 2}$)
(cf. \cite[Notation 3.4]{NW12}).

Now, let $\mathfrak{X}$ be a normal affine variety defined by
a ring $\mathcal{O}_{\mathfrak{X}}$, and suppose that 
a collection of morphisms 
$f_i:\mathfrak{X}\to \mathbf{P}^1-\{0,1,\infty\}$
and $c_i\in{\mathbb Z}$ $(i=1,\dots,n)$
satisfy a multi-linear relation (called the tensor criterion for
a functional equation)
\begin{equation}
\label{tensor-criterion}
\sum_{i=1}^n c_i \,
f_i^\ast(\widehat{\kappa_{\otimes m}}(\varphi))
\equiv
0
\qquad
\text{in} \ \ 
{\overline{\mathcal{O}}_{\mathfrak{X}}^{\times}}^{\otimes n-2} 
\otimes 
\bigl(\overline{\mathcal{O}}_{\mathfrak{X}}^{\times} \wedge 
\overline{\mathcal{O}}_{\mathfrak{X}}^{\times} 
\bigr),
\end{equation}
where $f_i^\ast: \mathcal{O}\to \mathcal{O}_{\mathfrak{X}}$ is 
the pull-back of functions on $U$ and 
$\overline{\mathcal{O}}_{\mathfrak{X}}^\times:=
\mathcal{O}_{\mathfrak{X}}^\times/{\mathbb C}^\times$.
For each $i=1,\dots,n$ and a topological 
path $\gamma:v{\leadsto} \xi$ on $\mathfrak{X}({\mathbb C})$,
the image $f_i(\gamma)$ forms a
path from $f_i(v)$ to $f_i(\xi)$
on $\mathbf{P}^1({\mathbb C})-\{0,1,\infty\}$.
Let  $\Lambda_{f_i(\gamma)}$ (or written also $\Lambda(f_i(\gamma))$) denote 
Chen's transport formal series (discussed in the proof of Lemma \ref{KeyId}) 
whose coefficients 
are iterated integrals along the path $f_i(\gamma)$. 
This is known as a group-like element of 
${\mathbb C}\langle\!\langle X,Y\rangle\!\rangle$ so that 
$\log(\Lambda_{f_i(\gamma)}^{-1})$
lies in the space of Lie formal series 
$\mathrm{Lie}_{\mathbb C}\langle\!\langle X,Y\rangle\!\rangle$ inside ${\mathbb C}\langle\!\langle X,Y\rangle\!\rangle$.
The homogeneous degree $m$ part $\mathrm{Lie}_m({\mathbb C})$
of $\mathrm{Lie}_{\mathbb C}\langle\!\langle X,Y\rangle\!\rangle$ is naturally
isomorphic to $\mathrm{gr}^m(\pi)\otimes {\mathbb C}$.
Define $\mathcal{L}_{\mathbb C}^\varphi(f_i(\xi),f_i(v);f_i(\gamma))$
(called the complex iterated integrals in \cite[Definition 4.4]{NW12})
to be the image of $\log(\Lambda_{f_i(\gamma)}^{-1})$
under the composition
\begin{equation}
\mathrm{Lie}_{\mathbb C}\langle\!\langle X,Y\rangle\!\rangle\twoheadrightarrow\mathrm{Lie}_m({\mathbb C})
\,{\overset \sim \to  }\, \mathrm{gr}^m(\pi)\otimes{\mathbb C}
\overset{\varphi_{\mathbb C}}\longrightarrow {\mathbb C}
\end{equation}
induced from $\varphi_{\mathbb C}:=\varphi\otimes{\mathbb C}$.
Then it follows from \cite[Theorem 4.13]{NW12}
that
\begin{equation}
\label{complex-FE}
\sum_{i=1}^n c_i\, \mathcal{L}_{\mathbb C}^\varphi(f_i(\xi),f_i(v);f_i(\gamma))=0.
\end{equation}
Next, to exhibit the $\ell$-adic Galois version,
suppose that $\mathfrak{X}$ and points $v$, $\xi$ on it together with
$f_i:\mathfrak{X}\to U$ $(i=1,\dots,n)$ are defined over  
a subfield $K\subset {\mathbb C}$, and assume that the tensor
condition (\ref{tensor-criterion}) is satisfied by them.
We fix a path system $\delta_i:\overrightarrow{01}{\leadsto} f_i(v)$
$(i=1,\dots,n)$.
For each pro-$\ell$ path $\gamma:v{\leadsto} \xi$ on $\mathfrak{X}_{\overline{K}}$
and $\sigma\in G_K$, we define
$\mathcal{L}^{\varphi(f_i),\vec{x}}_{\rm nv }\bigl(f_i(\xi),f_i(v);f_i(\gamma)\bigr)(\sigma)$
(which is called the naive $\ell$-adic iterated integral in \cite[Definition 4.7]{NW12})
to be the image of the Lie formal series
$\log\bigl(
(\delta_i\cdot
f_i(\gamma)\cdot \sigma(f_i(\gamma)^{-1})
\cdot \delta_i^{-1})^{-1}
\bigr)$
under the composition
\begin{equation}
\mathrm{Lie}_{{\mathbb Q}_\ell}
\langle\!\langle X,Y\rangle\!\rangle\twoheadrightarrow\mathrm{Lie}_m({\mathbb Q}_\ell)
\,{\overset \sim \to  }\, \mathrm{gr}^m(\pi)\otimes{\mathbb Q}_\ell\overset{\varphi_{{\mathbb Q}_\ell}}
\longrightarrow {\mathbb Q}_\ell
\end{equation}
induced from $\varphi_{{\mathbb Q}_\ell}:=\varphi\otimes{\mathbb Q}_\ell$.
  (The prefix $\vec{x}$ of
$\mathcal{L}^{\varphi(f_i),\vec{x}}_{\rm nv }$
is to designate dependency on $\vec{x}=(l_0,l_1)$,
the initially fixed free
generator of $\pi$; note that $X=\log(l_0), Y=\log(l_1)$
in the $\ell$-adic Galois case).
Then, \cite[Theorem 4.14]{NW12} implies
\begin{equation}
\label{ladic-FE}
\sum_{i=1}^n
c_i\,
\mathcal{L}^{\varphi(f_i),\vec{x}}_{\rm nv }\bigl(
f_i(\xi),f_i(v);f_i(\gamma)
\bigr)(\sigma)
=E(\sigma,\gamma)
\qquad (\sigma\in G_K),
\end{equation}
where the error term $E(\sigma,\gamma)$ is a function
of $\sigma\in G_K$ and of $\gamma:v{\leadsto} \xi$
satisfying a certain condition on small variation over 
$(\sigma,\gamma)$ (\cite[p.276]{NW12}).

The map $\varphi_3$ employed in 
(\ref{complex-basic}), (\ref{ladic-basic})
corresponds to the trilogaritm, and
is the special case $m=3$ of $\varphi_m:\mathrm{gr}^m(\pi)\to{\mathbb Z}$
defined as
the dual element to $\mathrm{ad}_X^{m-1}(Y)
=[X[...[X[X,Y]]..]]$ with respect to 
the Hall basis (for the order $X<Y$) 
of the free Lie algebra generated by $X,Y$. 
To obtain a polylogarithmic identity in the case $\varphi=\varphi_m$,
we first translate the ``iterated integrals along $f_i(v){\leadsto} f_i(\xi)$''
appearing in the left hand sides of (\ref{complex-FE}) and (\ref{ladic-FE})
into the terms of those along $\overrightarrow{01}{\leadsto} f_i(v)$ 
and $\overrightarrow{01}{\leadsto} f_i(\xi)$ under the natural network of paths
composed of $\{f_i(\gamma)\}_{i=1}^n \cup\{\delta_i\cdot f_i(\gamma)\}_{i=1}^n$.
This is achieved by what is called a ``polylog BCH formula'' elaborated in 
\cite[Proposition 5.9]{NW12}.
The second task is to evaluate the error term $E(\sigma,\gamma)$ 
of the right hand side of (\ref{ladic-FE}), which is figured out in 
\cite[Corollary 5.8]{NW12}.
We refer the reader to \cite[Sect.5]{NW12} for more details 
of these procedures when $\varphi$ is of the form $\varphi_m$ $(m\in{\mathbb N}_{\ge 2})$.

On the other hand, when $\varphi$ is chosen to be a general character looking at non-polylogarithmic 
(or to say, multiple polylogarithmic) coefficients, then there occur more
complicated computational procedures. We expect future studies for them.

Finally, we would like to present the following list of typos in the previous papers 
\cite{NW12}, \cite{NW20}, many pieces of which have been found during
the course of our present collaboration.
Since these papers are not only crucial but also indispensable to test computations 
of our above investigation, we hope to be allowed to append the list here below:

\smallskip\noindent
$\blacktriangleright$
{\tt Misprints in \cite{NW12}:}
\\
p.284, line $-11$: The LHS of the displayed formula should read:
$
%%Correction!!%%%%%%%%%%%%%%%%%%%
1-\frac{1}{2}(e^{(\log z) X}-1)%^2
+\frac{1}{3}(e^{(\log z)X}-1)^2%^3
%%%%%%%%%%%%%%%%%%%%%%%%%%%%%%%%%%%%%%%
-+\cdots
$
\\
p.286, line 6: The two Galois groups $G_\ast$
should respectively have subscripts
$\ast={K(\mu_{\ell^\infty},z^{1/\ell^{\infty}})}$ and
$\ast={K(\mu_{\ell^\infty},z^{1/\ell^{\infty}},(1-z)^{1/\ell^{\infty}})}$.
\\
p.287, in (iii$)_{\mathbb C}$: $f_i(x)$ should read $f_i(v)$, and
$x{\leadsto} z$ should read $v{\leadsto} z$.
\\
p.288, in the last line of (iii$)_\ell$, 
$p:x{\leadsto} z$ should read $p:v{\leadsto} z$
\\
p.294: In the formula (5.13), the denominator of RHS should read:
$
\prod_{a=0}^{\ell^n-1}
(1-\zeta_{\ell^n}^{a})^{\frac{a^{2k-1}}{\ell^n}}.
$
\\
p.300: The last line of calculation of $E(\sigma,\gamma)$ 
(displayed in the middle of page) should read:
$$
=\varphi_{2,\vec{x}}
\left(-\frac{1}{2}\rho_{1-z}(\sigma)
[\log\,x,\log\,y])\right)
=-\frac{1}{2}\rho_{1-z}(\sigma).
$$ 
({\it Note}: The subsequent displayed formula reflects 
$(6.20)$ with both sides multiplied by $-1$ so that
the formula (6.22) is itself correct.) 
\\
p.300, line $-1$: $f_1(v)=\overrightarrow{\infty 1}$ should read $f_2(v)=\overrightarrow{\infty 1}$.
\\
p.303, line $-7$: The latter equality should read 
$\delta f_2(y)\delta^{-1}=y$.
\\
p.304, line 2: The RHS of the displayed identity should read
$
=-\left(\frac{-\chi }{e^{-\chi t}-1}\right)
-\left(\frac{e^{{{\mathsf{L}}_0}t}}{e^{-t}-1}\right).
$
\\
p.305, line 4: The second identity should read
$f_2(x_{14})=x_{13}=(yx)^{-1}$.
\\
({\it Note}: The forgetful map $f_2:M_{0,5}\to M_{0,4}=
\mathbf{P}^1-\{0,1,\infty\}$ sends braid-like generators
in such a way that
$x_{ij}$ ($1\le i<j\le 5)$ of $\pi_1(M_{0,5}({\mathbb C}),\overrightarrow{v})$
is mapped to $x_{kl}$ ($1\le k<l\le 4$) of $\pi_1(M_{0,4}({\mathbb C}),\overrightarrow{v})$ by $k:=i-1$, $l:=j-1$,
except the cases with $f_2(x_{ij})=1$ ($2\in \{i,j\}$)
or with
$f_2(x_{1j})=x_{1,j-1}$ $(j=3,4,5)$. We have
$x_{12}x_{13}x_{23}=1$ in $\pi_1(M_{0,4}({\mathbb C}),\overrightarrow{v})$
freely generated by $x:=x_{12}$ and $y:=x_{23}$.)
\\
p.305, line 9: The first line of the table should read:
$$
%\begin{array}{c|ccccc}
\underline{
\quad \sharp \quad
\mid 
\mathrm{gr}^2_\Gamma f_{5\ast}(\sharp) 
%& 
\quad-\mathrm{gr}^2_\Gamma f_{4\ast}(\sharp) 
\quad %& 
\mathrm{gr}^2_\Gamma f_{3\ast}(\sharp) 
%&
\quad-\mathrm{gr}^2_\Gamma f_{2\ast}(\sharp) 
\quad %& 
\mathrm{gr}^2_\Gamma f_{1\ast}(\sharp) 
\ 
}
%\\
%\hline
%\end{array}
$$

\smallskip\noindent
$\blacktriangleright$
{\tt Misprints in \cite{NW20}:}
\\
p.603, line $-11$: ``$(-1)^{m-1}$-multiple'' should read ``$(-1)^{m}$-multiple''.
\\
p.603, line $-10$: The LHS of $(**)$ should read 
$(-1)^{m} {\ell}i_m(z,\gamma)$.
\\
p.603, line $-8$: The equality should read 
$d_{i+1}={\tilde{\chi}_{i+1}^z(\sigma)}{/i!}$ ($i\ge 0$).

\end{document}